\documentclass{svjour3}
\usepackage{amssymb}

\usepackage{epic}
\usepackage{amsmath}
\usepackage{amsfonts}
\usepackage{pstricks}

\usepackage{float}
\usepackage{subfig}
\usepackage{graphicx}
\usepackage{multirow}
\usepackage{adjustbox}
\usepackage{color}
\usepackage[subfigure]{tocloft}
\usepackage{subfig}
\usepackage[utf8]{inputenc}
\usepackage{soul}
\usepackage[hidelinks]{hyperref}

\allowdisplaybreaks
\numberwithin{equation}{section}
\usepackage[normalem]{ulem}
\newenvironment{rcases*}{\left.\begin{array}{ll}}{\end{array}\right\rbrace}

\def\bigO{\mathcal{O}}
\def\mig{1/2}

\def\bigO{\mathcal{O}}
\journalname{Journal of Scientific Computing}

\begin{document}
\title{Central WENO schemes through a global average weight}
\titlerunning{Central WENO schemes through a global average weight}
\author{Antonio~Baeza \and Raimund~B\"{u}rger \and Pep~Mulet \and David~Zor\'{\i}o}
\authorrunning{A.\ Baeza, R.\ B\"{u}rger,  P.\ Mulet and D.\ Zor\'{\i}o}
\date{}

\institute{A.\ Baeza \and P.\ Mulet  \at 
  Departament de Matem\`{a}tiques \\  Universitat de
  Val\`{e}ncia \\ E-46100 Burjassot,      Spain \\ 
  \email{antonio.baeza@uv.es, mulet@uv.es} 
  \and
  R.\ B\"{u}rger \at  
  CI$^2$MA \& Departamento de Ingenier\'{\i}a Matem\'{a}tica \\ 
   Universidad de Concepci\'{o}n  \\ 
   Casilla 160-C,  Concepci\'{o}n,    Chile\\ 
   \email{rburger@ing-mat.udec.cl} 
   \and 
   D.\ Zor\'{\i}o \at 
    CI$^2$MA,  
   Universidad de Concepci\'{o}n  \\ 
   Casilla 160-C,  Concepci\'{o}n,    Chile\\ 
   \email{dzorio@ci2ma.udec.cl}  
}

\thispagestyle{empty}

	\noindent This version of the article has been accepted for publication, after a peer-review process, and is subject to Springer Nature’s AM terms of use, but is not the Version of Record and does not reflect post-acceptance improvements, or any corrections. The Version of Record is available online at:

\noindent \url{https://doi.org/10.1007/s10915-018-0773-z}

\newpage
\setcounter{page}{1}
\maketitle

\begin{abstract}
  A novel central weighted essentially non-oscillatory (central WENO; CWENO)-type scheme for the construction of high-resolution approximations to discontinuous 
   solutions to hyperbolic systems of conservation laws 
  is presented. This procedure is based on the construction of a global average weight using  the whole set of Jiang-Shu smoothness indicators associated to every candidate stencil. By this device one does not to have to rely on ideal weights, which, under certain stencil arrangements and interpolating point locations, do not define a convex combination of the lower-degree interpolating polynomials of the corresponding sub-stencils. Moreover, this procedure also prevents some cases of accuracy loss near smooth extrema that are experienced by classical WENO and CWENO schemes. These properties result in  a more flexible scheme that overcomes these issues, at the cost of only a few additional computations 
  with respect to classical WENO schemes and with a smaller cost than classical CWENO schemes.  Numerical examples illustrate that the proposed CWENO schemes outperform both the traditional WENO and the original CWENO schemes. 

  \keywords{
    Finite difference schemes, central WENO schemes, global average weight.
  }
\end{abstract}

\section{Introduction}

\subsection{Scope} 

Weighted Essentially Non-Oscillatory (WENO) schemes
\cite{JiangShu96,LiuOsherChan94} have been  widely used in the
literature, especially in the context of  the approximation of  discontinuous solutions 
 to hyperbolic systems of 
conservation laws. The main feature of
the WENO  procedure is based on the fact that a reconstruction polynomial
can be decomposed as a certain convex combination of reconstruction
polynomials of lower order, provided they are evaluated at points within
a certain range. This property is attained in the case of the
well-known classical odd-order WENO schemes, both when the
 WENO procedure is applied for the interpolation of a function from point values 
  (as is expounded in \cite[Sect.~2.1]{shu09}) and for the reconstruction of a function  from cell averages 
   (see \cite[Sect.~2.2]{shu09}).
    The latter usage is more relevant to the numerical solution of 
    conservation laws. The weights used to ponder the contribution of each lower-order polynomial
    depend on the interpolating point and are known as ideal weights. WENO schemes define nonlinear 
    weights based on the ideal weights so as to construct an essentially non-oscillatory interpolant.
However, there are some circumstances in which the ideal linear weights are negative and thus the non-linear WENO weights do not satisfy the required properties, namely to attain the optimal order under smoothness assumptions and  to be essentially non-oscillatory when a discontinuity crosses the stencil. This is a well-known problem and strategies to solve it are summarized in \cite[Sect.~2.3.3]{shu09}.  
Therefore, since there are some practical situations in which classical WENO schemes are not suitable for use, some authors have proposed   solutions to overcome the aforementioned issues. In particular, Levy, Puppo and Russo introduced a central WENO (CWENO) scheme in \cite{LevyPuppoRusso}. We 
	will refer to their method as CWENO-LPR scheme.

In this work we propose a central WENO 
	scheme constructed following a different approach, based on a global average weight which does not depend on the ideal weights and is built using only the classical Jiang-Shu smoothness
indicators \cite{JiangShu96} that would be considered to compute the classical WENO
weights. Hence, it suffices to consider only
two additional items: on the one hand, the global average weight,
which is defined using the smoothness indicators through elementary
operations; and on the other the evaluation of the reconstruction
polynomial from the whole stencil.
Therefore, in terms of computational cost the new scheme is slightly more expensive than Jiang-Shu's WENO scheme, but much cheaper
	than the CWENO-LPR scheme of \cite{LevyPuppoRusso}, as the latter involves the computation of
	an additional  global smoothness indicator.

Finally, along the paper it will be shown that this procedure is also capable to overcome in some cases
 the issue of loss of  order of accuracy  near smooth extrema associated with  the original WENO and CWENO schemes.

\subsection{Related work} 

To further put the paper into the proper perspective, we recall that WENO schemes build on the previously 
 introduced family of essentially non-oscillatory (ENO) schemes that are based on selecting the least oscillatory 
 polynomial for reconstruction (among several available candidates defined by their respective stencils), 
  see Harten et al.\ \cite{Harten1987} and Shu and Osher \cite{ShuOsher89,ShuOsher1989}. The underlying idea of WENO schemes, namely to utilize a weighted 
   combination of these polynomials, was introduced in  \cite{LiuOsherChan94} and  put into a general framework 
    to construct arbitrary-order accurate finite difference schemes in~\cite{JiangShu96}. These schemes have 
   gained a vast amount of popularity and interest. For general information and references also to applications we refer to 
    review articles and handbook entries including \cite{shu98,shu09,zhangshu16}. 

The concept of central WENO (CWENO) schemes was advanced first by
Levy, Puppo and Russo 
in \cite{LevyPuppoRusso99} and later modified in \cite{LevyPuppoRusso} to define a compact CWENO scheme where the 
reconstruction polynomial is based on the information of the whole stencil, with the addition of an ideal weight associated to the reconstruction.
This modification allows to attain the optimal order for any convex combination of such weights, yielding a much more versatile scheme. See also \cite{Capdeville2008,CraveroSemplice} and references therein for further details regarding the aforementioned schemes.

\subsection{Outline of the paper} 

The remainder of the  paper is divided as follows: in Section~\ref{enm}, we briefly present the context in which we will stress out the performance of the proposed scheme along the paper. Section~\ref{cweno} is devoted to the description of the novel CWENO scheme in full detail. A motivation for developing
the new method is presented in Section~\ref{mot}, where we show some cases in which classical WENO schemes fail to provide a satisfactory strategy to perform spatial reconstructions. Section~\ref{nf} is focused on the formulation of our new scheme. Finally, in Section~\ref{aa} some theoretical results involving the accuracy of the weights and the reconstructions through our scheme are shown. Next, in Section \ref{ne}, several numerical tests are  presented in order to validate with numerical evidence the theoretical considerations drawn in the previous sections regarding the scheme presented 
in this paper. On one hand, Section~\ref{acte} is devoted to an extensive accuracy analysis; on the other hand, Section~\ref{clex} is focused on  several tests   to check the behaviour of the proposed scheme in shock problems from hyperbolic conservation laws,  and to compare them with the results obtained through the classical WENO and CWENO schemes. Finally, in Section \ref{cnc} some conclusions are drawn.

\section{Equations and numerical method}\label{enm}

Although WENO reconstructions are   not directly related to numerical schemes for a specific type of PDE, 
  we  focus on hyperbolic conservation laws. Therefore, we will briefly describe in this section the equations and their discretization procedure.
We consider hyperbolic systems of~$\nu$ scalar 
conservation laws in $d$~space dimensions: 
\begin{align}\label{hcl}
    \boldsymbol{u}_t+\sum_{i=1}^d \boldsymbol{f}^i(\boldsymbol{u})_{x_i}&= \boldsymbol{0},
     \quad (\boldsymbol{x},t)\in\Omega\times\mathbb{R}^+\subseteq\mathbb{R}^d\times\mathbb{R}^+, 
      \quad \boldsymbol{x} = (x_1,\ldots,x_d), 
\end{align} 
where $\boldsymbol{u}= \boldsymbol{u} ( \boldsymbol{x},t)\in\mathbb{R}^{\nu}$ is the sought solution, 
 $\boldsymbol{f}^i: \mathbb{R}^{\nu} \rightarrow\mathbb{R}^{\nu}$ are given flux density vectors, and 
\begin{align*}  
 \boldsymbol{u}=\begin{pmatrix} 
    u_1 \\
    \vdots \\
    u_{\nu}
    \end{pmatrix},\quad \boldsymbol{f}^i= \begin{pmatrix}       f^i_1 \\
      \vdots \\
      f^i_{\nu}
      \end{pmatrix}, \quad i=1, \dots, d;  \quad 
     \boldsymbol{f}=\begin{bmatrix} \boldsymbol{f}^1&\dots& \boldsymbol{f}^d
      \end{bmatrix}.
     \end{align*} 
System \eqref{hcl} is complemented with  the initial condition 
\begin{align*}
u(\boldsymbol{x},0)=\boldsymbol{u}_0(\boldsymbol{x}),  \quad  \boldsymbol{x}\in\Omega, 
\end{align*} 
and prescribed boundary conditions. 

To describe the spatial discretization, we introduce a  Cartesian grid $\mathcal{G}$ formed by points (cell centers) $\smash{\boldsymbol{x}=\boldsymbol{x}_{j_1,\dots,j_d}=((j_1-\frac{1}{2})h,\dots,(j_d-\frac{1}{2})h)\in\mathcal{G}}$ for $h>0$. 
 In what follows, we use the index vector $\boldsymbol{j} = (j_1, \dots, j_d)$,  let $\boldsymbol{e}_i$ denote the $i$-th $d$-dimensional unit vector, and assume that $J$ is the set of all indices $\boldsymbol{j}$ for which point values need to be updated. We then define 
 \begin{align*} 
  \boldsymbol{U} (t) := \bigl( \boldsymbol{u} ( \boldsymbol{x}_{\boldsymbol{j}}, t ) \bigr)_{\boldsymbol{j}\in J}. 
  \end{align*}   
To solve \eqref{hcl}
we utilize  the Shu-Osher finite difference scheme \cite{ShuOsher1989} with upwind spatial reconstructions of the flux function that are incorporated  into numerical flux vectors  $\smash{\boldsymbol{\hat f}}^{i}$ through a Donat-Marquina flux-splitting \cite{DonatMarquina96}. Thus, the contribution to the flux divergence in the coordinate~$x_i$ at point
 $\smash{\boldsymbol{x}=\boldsymbol{x}_{\boldsymbol{j}}}$  is given by 
\begin{align*} 
\boldsymbol{f}^i( \boldsymbol{U})_{x_i}(\boldsymbol{x}_{\boldsymbol{j}}, t) &  \approx\frac{1}{h} \Bigl( 
  \smash{\boldsymbol{\hat f}}^i_{\boldsymbol{j} + \frac{1}{2} \boldsymbol{e}_i} \bigl(\boldsymbol{U} ( t)\bigr)-
   \smash{\boldsymbol{\hat f}}^i_{ \boldsymbol{j} - \frac{1}{2} \boldsymbol{e}_i} \bigl( \boldsymbol{U}( t) \bigr) \Bigr). 
 \end{align*}
As a particular case of interest we consider WENO reconstructions \cite{JiangShu96} of order $2r+1$. 
To specify the time discretization, we write the semi-discrete scheme compactly as 
\begin{align*} 
\frac{\mathrm{d}} {\mathrm{d} t} \boldsymbol{U} (t) = \boldsymbol{\mathcal{L}} ( \boldsymbol{U} (t) ), 
\quad  \boldsymbol{\mathcal{L}} \bigl( \boldsymbol{U} (t) \bigr) =  \bigl( \mathcal{L}_{\boldsymbol{j}}  ( \boldsymbol{U} (t) ) \bigr)_{\boldsymbol{j} 
 \in J}, 
\end{align*} 
where we define 
\begin{align*} 
\mathcal{L}_{\boldsymbol{j}}  ( \boldsymbol{U} (t) ):= \frac{1}{h} \sum_{i=1}^d \Bigl( 
  \smash{\boldsymbol{\hat f}}^i_{\boldsymbol{j} + \frac{1}{2} \boldsymbol{e}_i} \bigl(\boldsymbol{U} ( t)\bigr)-
   \smash{\boldsymbol{\hat f}}^i_{ \boldsymbol{j} - \frac{1}{2} \boldsymbol{e}_i} \bigl( \boldsymbol{U}( t) \bigr) \Bigr)
\end{align*} 
(with suitable modifications for boundary points). 

For the time discretization, we use the third-order  Runge-Kutta TVD scheme proposed in \cite{ShuOsher89}. 
Assume that $\boldsymbol{U}^n := \boldsymbol{U} ( t^n)$ 
    is given and $\boldsymbol{U}^{n+1} = \boldsymbol{U} ( t^{n+1})$ is sought, 
    where $t^{n+1} = t^n + \Delta t$. Then  this scheme is defined as follows:
  \begin{align*}
 \begin{split} 
      \boldsymbol{U}^{(1)}&= \boldsymbol{U}^n-\Delta t \boldsymbol{\mathcal{L}}( \boldsymbol{U}^n),\\
      \boldsymbol{U}^{(2)}&= \frac{3}{4} \boldsymbol{U}^n+\frac{1}{4} \boldsymbol{U}^{(1)}-\frac{1}{4}\Delta
      t\boldsymbol{\mathcal{L}}(\boldsymbol{U}^{(1)}),\\
      \boldsymbol{U}^{n+1}&=\frac{1}{3}\boldsymbol{U}^n+\frac{2}{3}\boldsymbol{U}^{(2)}-\frac{2}{3}\Delta
      t\boldsymbol{\mathcal{L}}(\boldsymbol{U}^{(2)}).
      \end{split} 
      \end{align*}

\section{Central weighted essentially non-oscillatory (CWENO)  scheme}\label{cweno}

\subsection{Motivation}\label{mot}
To motivate the novel approach, which can be considered as an alternative CWENO scheme, let us focus on $d=1$ space dimension, 
and for ease of notation on a scalar equation ($\nu =1$). We drop
the $t$-dependence of $u$ for simplicity.

The key ingredient for obtaining highly accurate schemes for hyperbolic
conservation laws is the use of reconstructions that, given some contiguous
cell averages of an assumedly unknown function, produce precise local
evaluations. For classical finite volume schemes, the reconstructions
act on the evolved cell averages of the solution to precisely approximate
 the values of the solution at cell interfaces, whereas for  finite
difference schemes \cite{ShuOsher1989} the reconstructions are applied
to split fluxes to obtain at the end highly accurate approximations to
flux derivatives in conservative form.

In both cases the problem can be stated as the reconstruction of point values 
from cell averages. In the case of a finite-difference discretization on a 
uniform mesh 
the point values of the fluxes are assumed to be the (sliding) cell averages of a
certain function $\alpha$:
\begin{equation}\label{eq:t_shuosher}
  f(u(x)) = \frac{1}{h}\int_{x - h/2}^{x+h/2}\alpha \bigl(\xi\bigr)
  \, \mathrm{d}\xi, 
\end{equation}
so that the sought numerical fluxes $\hat{f}_{j+1/2}$ correspond to approximations
of the point values $\alpha(x+\frac h2)$
(see \cite{LiuOsherChan94}). Thus, we describe this problem from an
interpolatory point of view.

We perform a slight change of notation and assume that 
  ${f_{j-r},\ldots,f_{j+r}}$ are cell averages of a function $f(x)$
  associated to a stencil 
  of $2r+1$ points, such that 
  \begin{align*} 
   f_{j+l}=\frac{1}{h}\int_{x_{j+l-1/2}}^{x_{j+l+1/2}}f \bigl(x\bigr)
   \, \mathrm{d}x, \quad l=-r,\dots,r,
 \end{align*}
 and one wishes to obtain an approximation 
    \begin{align*} 
    \hat{f}_{j+\tau}=f \bigl(x_{j+\tau}\bigr)+\bigO(h^{2r+1}) \quad \text{for $0\leq\tau<1$}, 
    \end{align*}
    taking into account as well the  possible discontinuities in the data of the stencil. The case $\tau=\mig$ is well known and is handled properly by the traditional WENO schemes \cite{JiangShu96}, originally proposed in \cite{LiuOsherChan94}.

    This reconstruction method considers polynomials $p_{i,j}$ reconstructing point values from the cell average data $f_{j-r+i},\ldots,f_{j+i}$, $0\leq i\leq r$. Then, the Jiang-Shu smoothness indicators associated to these polynomials are computed as
    \begin{equation}\label{eq:si} 
      I_{i,j}:=\sum_{k=1}^r\int_{x_{j-1/2}}^{x_{j+1/2}}h^{2k-1} \bigl( p_{i,j}^{(k)}(x) \bigr)^2 \,\mathrm{d}x,\quad0\leq i\leq r.
    \end{equation}

    This scheme now uses the fact that there exists a convex combination $c_0,\ldots,c_r$, (namely, such that $c_0,\ldots,c_r>0$ and $\sum_{i=0}^{r+1}c_i=1$), called \textit{ideal linear weights}, such that $c_0p_{0,j}(x_{j+\mig})+\cdots+c_rp_{r,j}(x_{j+\mig})=p_j(x_{j+\mig})$, with $p_j$ the reconstruction polynomial associated to the whole stencil $\{f_{j-r},\ldots,f_{j+r}\}$. See for instance \cite[Proposition 2]{SINUM2011} for the explicit expression of $c_i$.

    Then, the \textit{non-linear weights}, which account for discontinuities in the data, are computed:
    $$\omega_{i,j}:=\frac{\alpha_{i,j}}{\alpha_{0,j}+\cdots+\alpha_{r,j}},\quad\alpha_{k,j}:=\frac{c_k}{(I_{k,j}+\varepsilon)^s},\quad0\leq i,k\leq r,\quad s>0,$$
    and the final WENO reconstruction is then given by
    $$q(x_{j+\mig})=\omega_{0,j}p_{0,j}(x_{j+\mig})+\cdots+\omega_{r,j}p_{r,j}(x_{j+\mig}).$$
 The general case of~$\tau$ is more complicated, since then  the ideal weights do not satisfy the same favourable properties as for  $\tau=\mig$. Let us analyze 
 as an example the particular case of   interest $r=2$ (namely, a fifth-order scheme). The three polynomials of degree~$2$ that interpolate three successive points of the stencil, and whose evaluations 
 at  $x_{j+\tau}=x_j+h\tau$  are 
 to be weighted within the WENO reconstruction, are given by  
\begin{align*}
  p_{0,j}(x_{j+\tau})&=\frac{-1+12\tau+12\tau^2}{24}f_{j-2}+\frac{1-24\tau-12\tau^2}{6}f_{j-1}+\frac{23+36\tau+12\tau^2}{24}f_j,\\
  p_{1,j}(x_{j+\tau})&=\frac{-1-12\tau+12\tau^2}{24}f_{j-1}+\frac{13-12\tau^2}{6}f_j+\frac{-1+12\tau+12\tau^2}{24}f_{j+1},\\
  p_{2,j}(x_{j+\tau})&=\frac{23-36\tau+12\tau^2}{24}f_j+\frac{1+24\tau-12\tau^2}{6}f_{j+1}+\frac{-1-12\tau+12\tau^2}{24}f_{j+2}.
\end{align*}
On the other hand, the result of interpolating on the whole stencil of five points and evaluating the resulting 
polynomial of degree 4 at~$x_{j+\tau}$ is 
\begin{align*}
  p_j(x_{j+\tau})&=\frac{9+200\tau-120\tau^2-160\tau^3+80\tau^4}{1920}f_{j-2} 
  \\ & \quad +\frac{-29-340\tau+360\tau^2+80\tau^3-80\tau^4}{480}f_{j-1}\\
  & \quad +\frac{1067-1320\tau^2+240\tau^4}{960}f_j\\
  & \quad +\frac{-29+340\tau+360\tau^2-80\tau^3-80\tau^4}{480}f_{j+1}\\
  &\quad +\frac{9-200\tau-120\tau^2+160\tau^3+80\tau^4}{1920}f_{j+2}. 
\end{align*}
The ideal weights  $c_0(\tau)$, $c_1(\tau)$ and $c_2(\tau)$ are rational expressions in~$\tau$ for  
 which 
\begin{align*} 
c_0(\tau)p_0(x_{j+\tau})+c_1(\tau)p_1(x_{j+\tau})+c_2(\tau)p_2(x_{j+\tau})=p(x_{j+\tau}).
\end{align*}
 In this case, we obtain 
\begin{align*}
  c_0(\tau)&=\frac{9+200\tau-120\tau^2-160\tau^3+80\tau^4}{-80+960\tau+960\tau^2},\\
  c_1(\tau)&=\frac{49-4548\tau^2+5360\tau^4-960\tau^6}{40-6720\tau+5760\tau^2},\\
  c_2(\tau)&=\frac{9-200\tau-120\tau^2+160\tau^3+80\tau^4}{-80-960\tau+960\tau^2}.
\end{align*}
Unfortunately, the ideal weights do not behave well for $0\leq\tau\leq1$, in the sense that they not only do  not satisfy the property $0\leq c_i(\tau)\leq1$, but also  are unbounded inside such range of $0\leq\tau<1$, which makes them unusable in practice. One can readily check that, for instance, $c_0$ has a pole at $\smash{\tau=-\frac{1}{2}+\frac{\sqrt{3}}{3}\approx0.07}$. However, we must point out that there are values $\tau$ which attain the desired properties involving $c_i(\tau)$, such as for $\tau=1/2$, corresponding to the ideal linear weights associated to the classical WENO schemes.

\subsection{Formulation}\label{cs}

The previous discussion related to the shortcomings of ideal weights motivates a new strategy to design weights for WENO reconstructions.  
This strategy is aimed to attain  the optimal order $2r+1$ when the stencil contains smooth data, and to reduce to $r$-th order when 
there is some avoidable discontinuity in the data.

\subsubsection{Classical Central WENO schemes (CWENO-LPR)}\label{ccws}

The  CWENO schemes introduced in \cite{LevyPuppoRusso} are described
		in this section. The scheme was originally described for the third-order case, and can be generalized 
		to arbitrary odd order $2r+1$. This approach uses the same smoothness indicators~$I_{i,j}$ as those defined 
in Equation \ref{eq:si}, then considers 
any  $r+2$ coefficients $c_0,\ldots,c_r,c_{r+1}$ in a convex
combination  and defines the following polynomial: 

$$p_{r+1,j}(x):=\frac{1}{c_{r+1}}\left(p_j(x)-\sum_{i=0}^rc_ip_{i,j}(x)\right),$$
with $p_j$ the reconstruction polynomial of the whole $(2r+1)$-point stencil.

Let us remark that there holds

$$\sum_{i=0}^{r+1}c_ip_{i,j}(x)=p_j(x),$$
and therefore $c_i$, $0\leq i\leq r+1$, act as ideal weights.

Then, the additional smoothness indicator associated to $p_{r+1,j}$ is computed

$$I_{r+1,j}:=\sum_{k=1}^{2r}\int_{x_{j-1/2}}^{x_{j+1/2}}h^{2k-1} \bigl( p_{r+1,j}^{(k)}(x) \bigr)^2 \,\mathrm{d}x.$$
The weights are then defined akin to those defined by Jiang-Shu, with the difference that now one additional polynomial, $p_{r+1,j}$, is included on the averaging:
$$\omega_{i,j}:=\frac{\alpha_{i,j}}{\alpha_{0,j}+\cdots+\alpha_{r,j}+\alpha_{r+1,j}},\quad\alpha_{k,j}:=\frac{c_k}{(I_{k,j}+\varepsilon)^s},\quad0\leq i,k\leq r+1,\quad s>0,$$
where $\varepsilon>0$ is a small number to  avoid  divisions by zero.

Then, the final reconstruction result at $x_{j+\tau}$ is defined as
$$\hat{f}_{j+\tau}=\sum_{i=0}^{r+1}\omega_{i,j}p_{i,j}(x_{j+\tau}).$$

\subsubsection{New formulation (CWENO)}\label{nf}

We next describe our proposal, whose main difference with respect to the classical approach is that
we now utilize instead the global average weight as defined in \cite{BaezaMuletZorio2016}, namely 
\begin{align} \label{omegajdef} 
\omega_j=\frac{(r+1)^2}{ \displaystyle \Biggl(\sum_{i=0}^r(I_i+\varepsilon)^{m}\Biggr)\Biggl(\sum_{i=0}^r\frac{1}{(I_i+\varepsilon)^m}\Biggr)},\quad m>0,
\end{align} 
where $\varepsilon>0$ is a small number to  avoid  divisions by zero. Moreover, loss of  accuracy at smooth extrema 
 is also avoided if one sets $\varepsilon=\bigO(h^2)$.
By  \cite[Prop.~2]{BaezaMuletZorio2016}, $\omega_j$ satisfies $0\leq\omega_j\leq1$ and, moreover, $\omega_j=1-\bigO(h^{2r})$ if the data from the stencil is smooth enough (assuming $\varepsilon=\bigO(h^2)$) and $\omega=\bigO(h^{2m})$ if there is a discontinuity.

Combining properly the weight~$\omega_j$ defined by \eqref{omegajdef}
with the $r+1$  polynomials of degree $r$, namely $p_{i,j}$ for $i=0,
\dots, r$, and the polynomial associated to the whole stencil of
degree $2r+1$, $p_j$, one can define the {\em reconstructed value} as
\begin{align*} 
\hat{f}_{j+\tau}=\omega_j  p_{j}(x_{j+\tau})+(1-\omega_j)q_j (x_{j+\tau}), \quad 
\text{where} \quad  q_j (x_{j+\tau}):=\sum_{i=0}^r\omega_{i,j} p_{i,j} (x_{j+\tau}), 
\end{align*} 
and the {\em subweights} are defined by 
\begin{align} \label{eq:subweights} 
 \omega_{i,j}:=\frac{\alpha_{i,j}}{\alpha_{0,j}+ \dots + \alpha_{r,j}},\quad\alpha_{k,j}:=\frac{c_{k}}{(I_{k,j}+\varepsilon)^s},\quad s>0,
 \end{align} 
where the constants  $c_k$ can be chosen such that $0<c_k<1$ and $c_0 + \dots + c_r=1$. This fact is the essential property of the so-called Central WENO  (CWENO) schemes described in Section \ref{ccws}, since, as pointed out in \cite{CraveroSemplice}, and extrapolating now the claim for reconstructions of arbitrary order, we impose the set of ideal weights in this case only by the condition $c_0 + \dots + c_r=1$, rather than the much more restrictive condition
$$\sum_{i=0}^rc_i p_{i,j} (x_{j+\mig})=p_{j}(x_{j+\mig}).$$
For instance, one can simply choose $c_k=\frac{1}{r+1}$ (arithmetic average) or the ideal weights of the case $\tau=\mig$ (the classical WENO schemes), which satisfy the required properties for any~$r$.  The above considerations imply  that 
\begin{align*} 
\hat{f}_{j+\tau}= f \bigl(x_{j+\tau} \bigr)+\bigO(h^{2r+1}) 
\end{align*} 
if there is smoothness and $\varepsilon=\bigO(h^2)$, and that 
\begin{align*} 
\hat{f}_{j+\tau}=f \bigl(x_{j+\tau} \bigr)+\bigO(h^{r+1})
\end{align*} 
otherwise.

Let us also point out that, unlike the classical procedure described in Section \ref{ccws}, this global average weight is independent of the remaining weights (the so-called {\it subweights}), and therefore it can be tuned at convenience through exponents without affecting the convexity properties of the remaining weights, as will be seen in Equation \eqref{remap} from Section \ref{aa}.

\subsection{Accuracy analysis}\label{aa}

We next study in full detail the accuracy of the reconstruction described
in Section~\ref{nf}  in terms of the choice of~$\varepsilon$ and the number of consecutive derivatives of~$f$  that vanish. This will be done 
 by  discussing   the parameters~$m$ and~$s$
that appear in \eqref{omegajdef} and \eqref{eq:subweights}, respectively. 
There are two plausible choices of~$\varepsilon$, namely $\varepsilon = \mathcal{O} ( h^2)$ (Choice~1) that avoids loss of accuracy at smooth extrema, 
 and  alternatively, $\varepsilon = \mathrm{const.}$  (Choice~2)  with an extremely small constant to  prevent divisions by zero, but that can 
  be   neglected in the accuracy analysis. In what follows we discuss the consequences for the accuracy analysis for either choice.

  For Choice~1 ($\varepsilon=\bigO(h^2)$), 
  the number of vanishing derivatives of~$f$ does not have any impact on $\omega_j$ or $\omega_{0,j}, \dots, \omega_{r,j}$  since these quantities  unconditionally satisfy
\begin{align*}
  \omega_j&=
  \begin{cases}
    1-\bigO(h^{2r}) & \textnormal{if the stencil is smooth,}\\
    \bigO(h^{2m}) & \textnormal{if a discontinuity crosses the stencil,}
  \end{cases} \\
  \omega_{i,j} &=
  \begin{cases}
    \bigO(1) & \textnormal{if substencil~$i$ is smooth, }\\
    \bigO(h^{2s}) & \textnormal{if a discontinuity crosses  substencil~$i$}
  \end{cases}
\end{align*}
(see \cite{SINUM2011} for further details). Therefore, since $\omega_{0,j}+ \dots + \omega_{r,j}=1$, one has
\begin{align} \label{qjeq} 
q_j(x_{j+\tau})=  f(x_{j+\tau})+
  \begin{cases}
   \bigO(h^{r+1}) & \textnormal{if the stencil is smooth,}\\
    \bigO(h^{\min\{2s,r+1\}}) & \textnormal{if a discontinuity crosses the stencil.}
  \end{cases}
\end{align}
Thus, to attain the optimal order~$r+1$ in the  second case of~\eqref{qjeq} it suffices to take
\begin{align} \label{sceil} 
s=\left\lceil (r+1)/2 \right\rceil.
\end{align} 
Finally, we obtain 
\begin{align*}
\hat{f}_{j+\tau} = f(x_{j+\tau})+ 
  \begin{cases}
   \bigO(h^{2r+1}) & \textnormal{if the stencil is smooth,}\\
   \bigO(h^{\min\{2m,r+1\}}) & \textnormal{if a discontinuity crosses the stencil.}
  \end{cases}
\end{align*} 
Thus, one can set
$m=\lceil (r+1)/2 \rceil$  to achieve the optimal order of accuracy, namely  
\begin{align} \label{eq:opt-ord} 
\hat{f}_{j+\tau} =  f(x_{j+\tau})+
  \begin{cases}
   \bigO(h^{2r+1}) & \textnormal{if the stencil is smooth,}\\
   \bigO(h^{r+1}) & \textnormal{if a discontinuity crosses the stencil.}
  \end{cases}
\end{align} 
This concludes the discussion of Choice~1 ($\varepsilon = \mathcal{O} ( h^2)$). 

In the alternative case of Choice~2 ($\varepsilon = \mathrm{const.} \ll 1$), 
 we must take into account the impact of smooth extrema in the accuracy of the weights. 
Then it is convenient to remap   the weight $\omega_j$ as
\begin{align}\label{remap} 
 \omega_j=(1-(1-\rho_j)^{s_1})^{s_2},
\quad \text{where we define} \quad  
\rho_j :=\frac{(r+1)^2}{\displaystyle \Biggl(\sum_{i=0}^r(I_i+\varepsilon)\Biggr) \Biggl(\sum_{i=0}^r\frac{1}{I_i+\varepsilon} \biggr)}.
\end{align} 
If $k_0$ is the maximum order of consecutive vanishing derivatives, that is 
 \begin{align*} 
 \frac{\mathrm{d}^k}{\mathrm{d} x^k} f(u)\bigr|_{x_{j+\tau}} =0, \quad 1 \leq k \leq k_0; \quad 
  \frac{\mathrm{d}^{k_0+1}}{\mathrm{d} x^{k_0+1}} f(u)\bigr|_{x_{j+\tau}} \neq 0,
 \end{align*}
then invoking once again  \cite[Prop.~2]{BaezaMuletZorio2016}, and defining $k_1:=\max\{r-k_0,0\}$, we get
\begin{align*}
  \omega_j =
  \begin{cases}
    1-\bigO(h^{2k_1s_1}) & \textnormal{if the stencil is smooth,}\\
    \bigO(h^{2s_2}) & \textnormal{if a discontinuity crosses the stencil.}
  \end{cases}
\end{align*}
Note that if $k_0\geq r$, then $k_1=0$, and $\omega_j$ does not  approximate~$1$,  so accuracy is unavoidably lost. We next show that when $k_0<r$ optimal accuracy (i.e.,  \eqref{eq:opt-ord})   
 can be recovered by setting properly~$s_1$ and~$s_2$.
Indeed, by choosing~$s$   by \eqref{sceil} (as for Choice~1) in the definition of  the weights $\omega_{i,j}$ 
 we obtain 
\begin{align*} 
q(x_{j+\tau})=f \bigl(x_{j+\tau} \bigr)+\bigO(h^{r+1}), 
\end{align*} 
whether the stencil is smooth or not.
 Taking into consideration the above remarks and assuming that there is smoothness, 
we deduce  that 
  \begin{align*}
    \hat{f}_{j+\tau}&=\omega_j p_j(x_{j+\tau})+(1-\omega_j)q_j (x_{j+\tau})\\
    &=(1-\omega_j) \bigl(f(x_{j+\tau})+\bigO(h^{2r+1}) \bigr)+\omega_j \bigl(f(x_{j+\tau})+\bigO(h^{r+1}) \bigr)\\
    &=f (x_{j+\tau} )+(1-\omega_j)\bigO(h^{2r+1})+\omega_j\bigO(h^{r+1})\\
    &=f(x_{j+\tau})+ \bigl(1-\bigO(h^{2k_1s_1})\bigr)\bigO(h^{2r+1})+\bigO(h^{2k_1s_1})\bigO(h^{r+1})\\
    &=f(x_{j+\tau})+\bigO(h^{2r+1})+\bigO(h^{2k_1s_1+r+1}).
  \end{align*}
  Thus, to achieve a $(2r+1)$-th accuracy order, one should impose
  $$2k_1s_1+r+1\geq 2r+1\Longleftrightarrow s_1\geq\frac{r}{2k_1}.$$
  Since $k_1 \geq 1$  whenever $k_0<r$,   the optimal exponent  is
  $s_1=\lceil r/2 \rceil$. 
 On the other hand, if  there is a discontinuity, the same reasoning as in the previous cases shows that the optimal parameter  is
  $s_2= \lceil (r+1)/2\rceil$. This ends the treatment of Choice~2. 
  
  We conclude the discussion on the order of accuracy by remarking that for either choice of~$\varepsilon$, both classical WENO and CWENO schemes of order $2r+1$ actually have order $r+1+|r-k|$, with $k=\min\{l\geq1\mid f^{(l)}(x_{j+\mig})\neq0\}-1$. Therefore, if we take into account our theoretical considerations, for $0<k<r$ the accuracy order of our scheme is the optimal, namely $2r+1$, whereas the accuracy order of the classical WENO and CWENO schemes in that case 
  is only $2r+1-k$.


\section{Numerical experiments}\label{ne}

The numerical experiments  are divided into two groups. In order to illustrate that the good behavior of our procedure is agnostic about the type of reconstructions to be weighted, we will use different type of reconstructions in each of these groups.

  The first group of numerical experiments (Examples~1 to~4) is devoted to algebraic problems where
   we test the performance of our CWENO method both on smooth problems and on
    discontinuous problems. The type of reconstructions used are 
    those  that interpret  the data of the  stencil as pointwise values of a certain 
    unknown function $f$. These experiments will be performed in Section~\ref{acte}.

    The second group of numerical experiments (Examples~5 to~7)  are shown in Section~\ref{clex} and involve numerical solutions of hyperbolic conservation laws through the Shu-Osher finite-difference method \cite{ShuOsher89,ShuOsher1989}. In this case the data of the stencils are the point values of the flux function, which, according to \eqref{eq:t_shuosher} are interpreted as the cell averages of some unknown function $\alpha$.
    
\subsection{Accuracy tests}\label{acte}

We now present some numerical tests, where we stress the performance of our scheme against that of  classical WENO schemes. These tests will be  focused on the quantitative behavior in presence of smooth extrema and discontinuities. To do so, we employ  an arbitrary precision library, MPFR \cite{MPFR}, using its C++ wrapper \cite{Holoborodko}, by setting a precision of $333$ bits ($\approx100$ digits).

\subsubsection*{Example 1: Smooth extrema analysis}

\begin{table}
\addtolength{\tabcolsep}{-4.5pt}
  \centering
  \begin{tabular}{|c|c|c|c|c|c|c|c|c|c|c|c|c|}
    \hline
    $r=1$ & \multicolumn{4}{c|}{WENO3} & \multicolumn{4}{c|}{CWENO3-LPR} & \multicolumn{4}{c|}{CWENO3}\\
    \hline
     & \multicolumn{2}{c|}{$k=0$} & \multicolumn{2}{c|}{$k=1$} & \multicolumn{2}{c|}{$k=0$} & \multicolumn{2}{c|}{$k=1$} & \multicolumn{2}{c|}{$k=0$} & \multicolumn{2}{c|}{$k=1$} \\
    \hline
    $n$ & Error & Ord. & Error & Ord. & Error & Ord. & Error & Ord. & Error & Ord. & Error & Ord. \\
    \hline
    40 & 8.97e-06 & --- & 1.56e-04 & --- & 1.02e-05 & --- & 1.56e-04 & --- & 2.86e-06 & --- & 1.56e-04 & --- \\
    \hline
    80 & 1.11e-06 & 3.01 & 3.91e-05 & 2.00 & 1.10e-06 & 3.20 & 3.91e-05 & 2.00 & 3.63e-07 & 2.98 & 3.91e-05 & 2.00 \\
    \hline
    160 & 1.38e-07 & 3.01 & 9.77e-06 & 2.00 & 1.27e-07 & 3.11 & 9.77e-06 & 2.00 & 4.56e-08 & 2.99 & 9.77e-06 & 2.00 \\
    \hline
    320 & 1.72e-08 & 3.00 & 2.44e-06 & 2.00 & 1.53e-08 & 3.06 & 2.44e-06 & 2.00 & 5.71e-09 & 3.00 & 2.44e-06 & 2.00 \\
    \hline
    640 & 2.15e-09 & 3.00 & 6.10e-07 & 2.00 & 1.87e-09 & 3.03 & 6.10e-07 & 2.00 & 7.15e-10 & 3.00 & 6.10e-07 & 2.00 \\
    \hline
    1280 & 2.68e-10 & 3.00 & 1.53e-07 & 2.00 & 2.31e-10 & 3.02 & 1.53e-07 & 2.00 & 8.94e-11 & 3.00 & 1.53e-07 & 2.00 \\
    \hline
    2560 & 3.35e-11 & 3.00 & 3.81e-08 & 2.00 & 2.87e-11 & 3.01 & 3.81e-08 & 2.00 & 1.12e-11 & 3.00 & 3.81e-08 & 2.00 \\
    \hline
    5120 & 4.19e-12 & 3.00 & 9.54e-09 & 2.00 & 3.58e-12 & 3.00 & 9.54e-09 & 2.00 & 1.40e-12 & 3.00 & 9.54e-09 & 2.00 \\
    \hline
  \end{tabular}
  \caption{Example 1: errors of schemes WENO3, CWENO3-LPR, and CWENO3.}
  \label{order_smooth_1}
\end{table}
\begin{table}
  \centering
  \begin{tabular}{|c|c|c|c|c|c|c|c|}
    \hline
    & $r=2$ & \multicolumn{2}{c|}{$k=0$} & \multicolumn{2}{c|}{$k=1$} & \multicolumn{2}{c|}{$k=2$} \\
    \cline{2-8}
    & $n$ & Error & Order & Error & Order & Error & Order \\
    \hline
    \multirow{11}{*}{\rotatebox[origin=c]{90}{WENO5}}
    & 40 & 7.35e-09 & --- & 1.44e-07 & --- & 1.90e-06 & --- \\
    \cline{2-8}
    & 80 & 2.31e-10 & 4.99 & 1.09e-08 & 3.72 & 1.98e-07 & 3.27 \\
    \cline{2-8}
    & 160 & 7.25e-12 & 5.00 & 7.44e-10 & 3.88 & 2.21e-08 & 3.16 \\
    \cline{2-8}
    & 320 & 2.27e-13 & 5.00 & 4.84e-11 & 3.94 & 2.60e-09 & 3.09 \\
    \cline{2-8}
    & 640 & 7.09e-15 & 5.00 & 3.08e-12 & 3.97 & 3.14e-10 & 3.05 \\
    \cline{2-8}
    & 1280 & 2.22e-16 & 5.00 & 1.95e-13 & 3.99 & 3.86e-11 & 3.02 \\
    \cline{2-8}
    & 2560 & 6.93e-18 & 5.00 & 1.22e-14 & 3.99 & 4.78e-12 & 3.01 \\
    \cline{2-8}
    & 5120 & 2.16e-19 & 5.00 & 7.66e-16 & 4.00 & 5.95e-13 & 3.01 \\
    \cline{2-8}
    & 10240 & 6.77e-21 & 5.00 & 4.79e-17 & 4.00 & 7.42e-14 & 3.00 \\
    \cline{2-8}
    & 20480 & 2.11e-22 & 5.00 & 3.00e-18 & 4.00 & 9.26e-15 & 3.00 \\
    \cline{2-8}
    & 40960 & 6.61e-24 & 5.00 & 1.87e-19 & 4.00 & 1.16e-15 & 3.00 \\
    \hline
    \multirow{11}{*}{\rotatebox[origin=c]{90}{CWENO5-LPR}}
    & 40 & 2.69e-08 & --- & 1.70e-06 & --- & 2.99e-07 & --- \\
    \cline{2-8}
    & 80 & 8.27e-10 & 5.02 & 7.46e-08 & 4.51 & 7.97e-08 & 1.91 \\
    \cline{2-8}
    & 160 & 2.57e-11 & 5.01 & 3.62e-09 & 4.37 & 1.26e-08 & 2.66 \\
    \cline{2-8}
    & 320 & 8.02e-13 & 5.00 & 1.93e-10 & 4.23 & 1.74e-09 & 2.86 \\
    \cline{2-8}
    & 640 & 2.50e-14 & 5.00 & 1.11e-11 & 4.13 & 2.28e-10 & 2.93 \\
    \cline{2-8}
    & 1280 & 7.81e-16 & 5.00 & 6.60e-13 & 4.07 & 2.91e-11 & 2.97 \\
    \cline{2-8}
    & 2560 & 2.44e-17 & 5.00 & 4.03e-14 & 4.03 & 3.68e-12 & 2.98 \\
    \cline{2-8}
    & 5120 & 7.63e-19 & 5.00 & 2.49e-15 & 4.02 & 4.63e-13 & 2.99 \\
    \cline{2-8}
    & 10240 & 2.38e-20 & 5.00 & 1.54e-16 & 4.01 & 5.80e-14 & 3.00 \\
    \cline{2-8}
    & 20480 & 7.45e-22 & 5.00 & 9.62e-18 & 4.00 & 7.26e-15 & 3.00 \\
    \cline{2-8}
    & 40960 & 2.33e-23 & 5.00 & 6.00e-19 & 4.00 & 9.08e-16 & 3.00 \\
    \hline
    \multirow{11}{*}{\rotatebox[origin=c]{90}{CWENO5}}
    & 40 & 5.65e-10 & --- & 6.16e-10 & --- & 1.42e-06 & --- \\
    \cline{2-8}
    & 80 & 1.78e-11 & 4.99 & 1.61e-11 & 5.26 & 1.53e-07 & 3.21 \\
    \cline{2-8}
    & 160 & 5.57e-13 & 5.00 & 1.29e-12 & 3.64 & 1.74e-08 & 3.14 \\
    \cline{2-8}
    & 320 & 1.74e-14 & 5.00 & 5.45e-14 & 4.57 & 2.06e-09 & 3.08 \\
    \cline{2-8}
    & 640 & 5.45e-16 & 5.00 & 1.94e-15 & 4.81 & 2.50e-10 & 3.04 \\
    \cline{2-8}
    & 1280 & 1.70e-17 & 5.00 & 6.43e-17 & 4.91 & 3.08e-11 & 3.02 \\
    \cline{2-8}
    & 2560 & 5.33e-19 & 5.00 & 2.07e-18 & 4.96 & 3.82e-12 & 3.01 \\
    \cline{2-8}
    & 5120 & 1.67e-20 & 5.00 & 6.57e-20 & 4.98 & 4.76e-13 & 3.01 \\
    \cline{2-8}
    & 10240 & 5.20e-22 & 5.00 & 2.07e-21 & 4.99 & 5.94e-14 & 3.00 \\
    \cline{2-8}
    & 20480 & 1.63e-23 & 5.00 & 6.48e-23 & 4.99 & 7.41e-15 & 3.00 \\
    \cline{2-8}
    & 40960 & 5.08e-25 & 5.00 & 2.03e-24 & 5.00 & 9.26e-16 & 3.00 \\
    \hline
  \end{tabular}
  \caption{Example~1: errors of schemes WENO5, CWENO5-LPR, and CWENO5.}
  \label{order_smooth_2}
\end{table}
\begin{table}
  \addtolength{\tabcolsep}{-3.5pt}
  \centering
  \begin{tabular}{|c|c|c|c|c|c|c|c|c|c|}
    \hline
    & $r=3$ & \multicolumn{2}{c|}{$k=0$} & \multicolumn{2}{c|}{$k=1$} & \multicolumn{2}{c|}{$k=2$} & \multicolumn{2}{c|}{$k=3$} \\
    \cline{2-10}
    & $n$ & Error & Order & Error & Order & Error & Order & Error & Order \\
    \hline
    \multirow{12}{*}{\rotatebox[origin=c]{90}{WENO7}}
    & 40 & 1.55e-12 & --- & 1.59e-10 & --- & 5.48e-09 & --- & 2.19e-07 & --- \\
    \cline{2-10}
    & 80 & 1.22e-14 & 6.99 & 2.12e-12 & 6.22 & 6.53e-11 & 6.39 & 1.37e-08 & 4.00 \\
    \cline{2-10}
    & 160 & 9.54e-17 & 7.00 & 3.04e-14 & 6.13 & 3.28e-13 & 7.64 & 8.54e-10 & 4.00 \\
    \cline{2-10}
    & 320 & 7.46e-19 & 7.00 & 4.52e-16 & 6.07 & 1.70e-14 & 4.27 & 5.34e-11 & 4.00 \\
    \cline{2-10}
    & 640 & 5.83e-21 & 7.00 & 6.89e-18 & 6.04 & 9.61e-16 & 4.15 & 3.34e-12 & 4.00 \\
    \cline{2-10}
    & 1280 & 4.55e-23 & 7.00 & 1.06e-19 & 6.02 & 3.68e-17 & 4.71 & 2.09e-13 & 4.00 \\
    \cline{2-10}
    & 2560 & 3.56e-25 & 7.00 & 1.65e-21 & 6.01 & 1.25e-18 & 4.87 & 1.30e-14 & 4.00 \\
    \cline{2-10}
    & 5120 & 2.78e-27 & 7.00 & 2.57e-23 & 6.00 & 4.09e-20 & 4.94 & 8.15e-16 & 4.00 \\
    \cline{2-10}
    & 10240 & 2.17e-29 & 7.00 & 4.01e-25 & 6.00 & 1.30e-21 & 4.97 & 5.09e-17 & 4.00 \\
    \cline{2-10}
    & 20480 & 1.70e-31 & 7.00 & 6.26e-27 & 6.00 & 4.11e-23 & 4.99 & 3.18e-18 & 4.00 \\
    \cline{2-10}
    & 40960 & 1.33e-33 & 7.00 & 9.78e-29 & 6.00 & 1.29e-24 & 4.99 & 1.99e-19 & 4.00 \\
    \cline{2-10}
    & 81920 & 1.04e-35 & 7.00 & 1.53e-30 & 6.00 & 4.04e-26 & 5.00 & 1.24e-20 & 4.00 \\
    \hline
    \multirow{12}{*}{\rotatebox[origin=c]{90}{CWENO7-LPR}}
    & 40 & 8.14e-12 & --- & 5.98e-09 & --- & 2.26e-07 & --- & 2.18e-07 & --- \\
    \cline{2-10}
    & 80 & 5.70e-14 & 7.16 & 9.50e-11 & 5.98 & 5.07e-09 & 5.48 & 1.36e-08 & 4.00 \\
    \cline{2-10}
    & 160 & 4.18e-16 & 7.09 & 1.50e-12 & 5.99 & 1.24e-10 & 5.36 & 8.51e-10 & 4.00 \\
    \cline{2-10}
    & 320 & 3.16e-18 & 7.05 & 2.35e-14 & 5.99 & 3.35e-12 & 5.21 & 5.32e-11 & 4.00 \\
    \cline{2-10}
    & 640 & 2.43e-20 & 7.02 & 3.68e-16 & 6.00 & 9.69e-14 & 5.11 & 3.33e-12 & 4.00 \\
    \cline{2-10}
    & 1280 & 1.88e-22 & 7.01 & 5.75e-18 & 6.00 & 2.91e-15 & 5.06 & 2.08e-13 & 4.00 \\
    \cline{2-10}
    & 2560 & 1.46e-24 & 7.01 & 8.99e-20 & 6.00 & 8.90e-17 & 5.03 & 1.30e-14 & 4.00 \\
    \cline{2-10}
    & 5120 & 1.14e-26 & 7.00 & 1.40e-21 & 6.00 & 2.75e-18 & 5.01 & 8.12e-16 & 4.00 \\
    \cline{2-10}
    & 10240 & 8.90e-29 & 7.00 & 2.20e-23 & 6.00 & 8.56e-20 & 5.01 & 5.07e-17 & 4.00 \\
    \cline{2-10}
    & 20480 & 6.95e-31 & 7.00 & 3.43e-25 & 6.00 & 2.67e-21 & 5.00 & 3.17e-18 & 4.00 \\
    \cline{2-10}
    & 40960 & 5.43e-33 & 7.00 & 5.36e-27 & 6.00 & 8.33e-23 & 5.00 & 1.98e-19 & 4.00 \\
    \cline{2-10}
    & 81920 & 4.24e-35 & 7.00 & 8.38e-29 & 6.00 & 2.60e-24 & 5.00 & 1.24e-20 & 4.00 \\
    \hline
    \multirow{12}{*}{\rotatebox[origin=c]{90}{CWENO7}}
    & 40 & 1.03e-13 & --- & 6.17e-13 & --- & 1.36e-11 & --- & 2.09e-07 & --- \\
    \cline{2-10}
    & 80 & 8.10e-16 & 6.99 & 4.85e-15 & 6.99 & 3.29e-14 & 8.69 & 1.31e-08 & 4.00 \\
    \cline{2-10}
    & 160 & 6.35e-18 & 7.00 & 3.81e-17 & 6.99 & 1.93e-16 & 7.41 & 8.19e-10 & 4.00 \\
    \cline{2-10}
    & 320 & 4.97e-20 & 7.00 & 2.98e-19 & 7.00 & 1.48e-18 & 7.03 & 5.12e-11 & 4.00 \\
    \cline{2-10}
    & 640 & 3.88e-22 & 7.00 & 2.33e-21 & 7.00 & 1.16e-20 & 6.99 & 3.20e-12 & 4.00 \\
    \cline{2-10}
    & 1280 & 3.03e-24 & 7.00 & 1.82e-23 & 7.00 & 9.09e-23 & 7.00 & 2.00e-13 & 4.00 \\
    \cline{2-10}
    & 2560 & 2.37e-26 & 7.00 & 1.42e-25 & 7.00 & 7.11e-25 & 7.00 & 1.25e-14 & 4.00 \\
    \cline{2-10}
    & 5120 & 1.85e-28 & 7.00 & 1.11e-27 & 7.00 & 5.56e-27 & 7.00 & 7.83e-16 & 4.00 \\
    \cline{2-10}
    & 10240 & 1.45e-30 & 7.00 & 8.68e-30 & 7.00 & 4.34e-29 & 7.00 & 4.89e-17 & 4.00 \\
    \cline{2-10}
    & 20480 & 1.13e-32 & 7.00 & 6.79e-32 & 7.00 & 3.39e-31 & 7.00 & 3.06e-18 & 4.00 \\
    \cline{2-10}
    & 40960 & 8.84e-35 & 7.00 & 5.30e-34 & 7.00 & 2.65e-33 & 7.00 & 1.91e-19 & 4.00 \\
    \cline{2-10}
    & 81920 & 6.90e-37 & 7.00 & 4.14e-36 & 7.00 & 2.07e-35 & 7.00 & 1.19e-20 & 4.00 \\
    \hline
  \end{tabular}
  \caption{Example~1: errors of schemes  WENO7, CWENO7-LPR, and CWENO7.}
  \label{order_smooth_3}
\end{table}
\begin{table}
  \addtolength{\tabcolsep}{-3.5pt}
  \centering
  \begin{tabular}{|c|c|c|c|c|c|c|c|c|c|c|c|}
    \hline
    & $r=4$ & \multicolumn{2}{c|}{$k=0$} & \multicolumn{2}{c|}{$k=1$} & \multicolumn{2}{c|}{$k=2$} & \multicolumn{2}{c|}{$k=3$} & \multicolumn{2}{c|}{$k=4$} \\
    \cline{2-12}
    & $n$ & Error & Ord. & Error & Ord. & Error & Ord. & Error & Ord. & Error & Ord. \\
    \hline
    \multirow{15}{*}{\rotatebox[origin=c]{90}{WENO9}}
    & 40 & 5.61e-16 & --- & 1.00e-13 & --- & 7.02e-12 & --- & 3.66e-10 & --- & 1.70e-09 & --- \\
    \cline{2-12}
    & 80 & 1.11e-18 & 8.99 & 3.26e-16 & 8.27 & 3.03e-14 & 7.86 & 8.24e-13 & 8.80 & 3.88e-11 & 5.45 \\
    \cline{2-12}
    & 160 & 2.17e-21 & 9.00 & 1.14e-18 & 8.16 & 1.68e-16 & 7.49 & 3.20e-15 & 8.01 & 2.64e-12 & 3.88 \\
    \cline{2-12}
    & 320 & 4.23e-24 & 9.00 & 4.19e-21 & 8.09 & 1.10e-18 & 7.26 & 4.95e-17 & 6.01 & 1.04e-13 & 4.67 \\
    \cline{2-12}
    & 640 & 8.27e-27 & 9.00 & 1.59e-23 & 8.05 & 7.83e-21 & 7.13 & 2.05e-19 & 7.91 & 3.58e-15 & 4.86 \\
    \cline{2-12}
    & 1280 & 1.62e-29 & 9.00 & 6.10e-26 & 8.02 & 5.86e-23 & 7.06 & 3.57e-21 & 5.85 & 1.17e-16 & 4.94 \\
    \cline{2-12}
    & 2560 & 3.16e-32 & 9.00 & 2.36e-28 & 8.01 & 4.48e-25 & 7.03 & 1.18e-22 & 4.92 & 3.73e-18 & 4.97 \\
    \cline{2-12}
    & 5120 & 6.16e-35 & 9.00 & 9.20e-31 & 8.01 & 3.46e-27 & 7.02 & 2.37e-24 & 5.64 & 1.18e-19 & 4.98 \\
    \cline{2-12}
    & 10240 & 1.20e-37 & 9.00 & 3.59e-33 & 8.00 & 2.69e-29 & 7.01 & 4.12e-26 & 5.84 & 3.70e-21 & 4.99 \\
    \cline{2-12}
    & 20480 & 2.35e-40 & 9.00 & 1.40e-35 & 8.00 & 2.09e-31 & 7.00 & 6.78e-28 & 5.93 & 1.16e-22 & 5.00 \\
    \cline{2-12}
    & 40960 & 4.59e-43 & 9.00 & 5.46e-38 & 8.00 & 1.63e-33 & 7.00 & 1.09e-29 & 5.96 & 3.63e-24 & 5.00 \\
    \cline{2-12}
    & 81920 & 8.97e-46 & 9.00 & 2.13e-40 & 8.00 & 1.28e-35 & 7.00 & 1.72e-31 & 5.98 & 1.14e-25 & 5.00 \\
    \cline{2-12}
    & 163840 & 1.75e-48 & 9.00 & 8.33e-43 & 8.00 & 9.96e-38 & 7.00 & 2.70e-33 & 5.99 & 3.55e-27 & 5.00 \\
    \cline{2-12}
    & 327680 & 3.42e-51 & 9.00 & 3.25e-45 & 8.00 & 7.78e-40 & 7.00 & 4.23e-35 & 6.00 & 1.11e-28 & 5.00 \\
    \cline{2-12}
    & 655360 & 6.68e-54 & 9.00 & 1.27e-47 & 8.00 & 6.08e-42 & 7.00 & 6.62e-37 & 6.00 & 3.47e-30 & 5.00 \\
    \hline
    \multirow{15}{*}{\rotatebox[origin=c]{90}{CWENO9-LPR}}
    & 40 & 1.35e-14 & --- & 4.28e-12 & --- & 3.51e-09 & --- & 3.52e-08 & --- & 1.10e-09 & --- \\
    \cline{2-12}
    & 80 & 2.61e-17 & 9.01 & 1.42e-14 & 8.24 & 2.88e-11 & 6.93 & 9.77e-11 & 8.49 & 1.24e-10 & 3.15 \\
    \cline{2-12}
    & 160 & 5.08e-20 & 9.01 & 5.02e-17 & 8.14 & 2.29e-13 & 6.97 & 3.41e-12 & 4.84 & 5.16e-12 & 4.59 \\
    \cline{2-12}
    & 320 & 9.91e-23 & 9.00 & 1.86e-19 & 8.08 & 1.80e-15 & 6.99 & 8.60e-14 & 5.31 & 1.80e-13 & 4.84 \\
    \cline{2-12}
    & 640 & 1.93e-25 & 9.00 & 7.06e-22 & 8.04 & 1.41e-17 & 7.00 & 1.56e-15 & 5.79 & 5.90e-15 & 4.93 \\
    \cline{2-12}
    & 1280 & 3.78e-28 & 9.00 & 2.72e-24 & 8.02 & 1.10e-19 & 7.00 & 2.58e-17 & 5.92 & 1.89e-16 & 4.97 \\
    \cline{2-12}
    & 2560 & 7.37e-31 & 9.00 & 1.05e-26 & 8.01 & 8.63e-22 & 7.00 & 4.13e-19 & 5.96 & 5.96e-18 & 4.98 \\
    \cline{2-12}
    & 5120 & 1.44e-33 & 9.00 & 4.11e-29 & 8.01 & 6.75e-24 & 7.00 & 6.54e-21 & 5.98 & 1.87e-19 & 4.99 \\
    \cline{2-12}
    & 10240 & 2.81e-36 & 9.00 & 1.60e-31 & 8.00 & 5.27e-26 & 7.00 & 1.03e-22 & 5.99 & 5.87e-21 & 5.00 \\
    \cline{2-12}
    & 20480 & 5.49e-39 & 9.00 & 6.25e-34 & 8.00 & 4.12e-28 & 7.00 & 1.61e-24 & 6.00 & 1.84e-22 & 5.00 \\
    \cline{2-12}
    & 40960 & 1.07e-41 & 9.00 & 2.44e-36 & 8.00 & 3.22e-30 & 7.00 & 2.52e-26 & 6.00 & 5.74e-24 & 5.00 \\
    \cline{2-12}
    & 81920 & 2.10e-44 & 9.00 & 9.53e-39 & 8.00 & 2.51e-32 & 7.00 & 3.94e-28 & 6.00 & 1.79e-25 & 5.00 \\
    \cline{2-12}
    & 163840 & 4.09e-47 & 9.00 & 3.72e-41 & 8.00 & 1.96e-34 & 7.00 & 6.16e-30 & 6.00 & 5.61e-27 & 5.00 \\
    \cline{2-12}
    & 327680 & 7.99e-50 & 9.00 & 1.45e-43 & 8.00 & 1.53e-36 & 7.00 & 9.63e-32 & 6.00 & 1.75e-28 & 5.00 \\
    \cline{2-12}
    & 655360 & 1.56e-52 & 9.00 & 5.68e-46 & 8.00 & 1.20e-38 & 7.00 & 1.50e-33 & 6.00 & 5.48e-30 & 5.00 \\
    \hline
    \multirow{15}{*}{\rotatebox[origin=c]{90}{CWENO9}}
    & 40 & 1.81e-17 & --- & 1.45e-16 & --- & 1.01e-15 & --- & 6.13e-12 & --- & 1.58e-09 & --- \\
    \cline{2-12}
    & 80 & 3.56e-20 & 8.99 & 2.85e-19 & 8.99 & 1.99e-18 & 8.99 & 1.07e-15 & 12.49 & 3.65e-11 & 5.44 \\
    \cline{2-12}
    & 160 & 6.97e-23 & 9.00 & 5.58e-22 & 8.99 & 3.90e-21 & 8.99 & 3.11e-19 & 11.74 & 2.50e-12 & 3.87 \\
    \cline{2-12}
    & 320 & 1.36e-25 & 9.00 & 1.09e-24 & 9.00 & 7.63e-24 & 9.00 & 3.36e-22 & 9.86 & 9.88e-14 & 4.66 \\
    \cline{2-12}
    & 640 & 2.67e-28 & 9.00 & 2.13e-27 & 9.00 & 1.49e-26 & 9.00 & 1.66e-25 & 10.98 & 3.41e-15 & 4.86 \\
    \cline{2-12}
    & 1280 & 5.21e-31 & 9.00 & 4.17e-30 & 9.00 & 2.92e-29 & 9.00 & 9.10e-29 & 10.83 & 1.11e-16 & 4.93 \\
    \cline{2-12}
    & 2560 & 1.02e-33 & 9.00 & 8.14e-33 & 9.00 & 5.70e-32 & 9.00 & 1.67e-31 & 9.09 & 3.56e-18 & 4.97 \\
    \cline{2-12}
    & 5120 & 1.99e-36 & 9.00 & 1.59e-35 & 9.00 & 1.11e-34 & 9.00 & 4.49e-34 & 8.54 & 1.12e-19 & 4.98 \\
    \cline{2-12}
    & 10240 & 3.88e-39 & 9.00 & 3.11e-38 & 9.00 & 2.17e-37 & 9.00 & 1.07e-36 & 8.72 & 3.53e-21 & 4.99 \\
    \cline{2-12}
    & 20480 & 7.58e-42 & 9.00 & 6.07e-41 & 9.00 & 4.25e-40 & 9.00 & 2.30e-39 & 8.85 & 1.11e-22 & 5.00 \\
    \cline{2-12}
    & 40960 & 1.48e-44 & 9.00 & 1.18e-43 & 9.00 & 8.29e-43 & 9.00 & 4.73e-42 & 8.93 & 3.46e-24 & 5.00 \\
    \cline{2-12}
    & 81920 & 2.89e-47 & 9.00 & 2.31e-46 & 9.00 & 1.62e-45 & 9.00 & 9.48e-45 & 8.96 & 1.08e-25 & 5.00 \\
    \cline{2-12}
    & 163840 & 5.65e-50 & 9.00 & 4.52e-49 & 9.00 & 3.16e-48 & 9.00 & 1.87e-47 & 8.98 & 3.38e-27 & 5.00 \\
    \cline{2-12}
    & 327680 & 1.10e-52 & 9.00 & 8.83e-52 & 9.00 & 6.18e-51 & 9.00 & 3.68e-50 & 8.99 & 1.06e-28 & 5.00 \\
    \cline{2-12}
    & 655360 & 2.16e-55 & 9.00 & 1.72e-54 & 9.00 & 1.21e-53 & 9.00 & 7.22e-53 & 9.00 & 3.30e-30 & 5.00 \\
    \hline
  \end{tabular}
  \caption{Example~1: errors of schemes WENO9, CWENO9-LPR, and CWENO9.}
  \label{order_smooth_4}
\end{table}

Let us consider the family of functions $f_k:\mathbb{R}\to\mathbb{R}$, $k\in\mathbb{N}$, given by
$$f_k(x)=x^{k+1}\mathrm{e}^x.$$
The function $f_k$ has a smooth extreme at $x=0$ of order $k$. We perform several tests involving different values of $k$ and $r$, where in each case the corresponding CWENO scheme with optimal parameters, $s_1=\lceil r / 2 \rceil$, $s_2=\lceil (r+1)/2 \rceil$, is compared against the classical WENO scheme and the CWENO-LPR scheme of the same order, with $\varepsilon=10^{-100}$ in all cases. We take the stencil
$x_j=(j-1/2)h$, $-r\leq j\leq r$, 
and perform the reconstruction at $x=0$ for different values of $h>0$. Since in this case one has a centered reconstruction point, we use the same ideal weights of the traditional WENO schemes to define the subweights. In this case, since the reconstructions are performed in the classical sense, namely, $\tau=\mig$, we choose as the non-linear subweights for the CWENO schemes those based on the ideal linear weights. We will perform accuracy tests for $1\leq r\leq 4$, namely, from (C)WENO3 to (C)WENO9, using the reconstructions from pointwise values to pointwise values, choosing in each case the smallest possible exponents~$s_1$ and~$s_2$ to attain the optimal order  \eqref{eq:opt-ord}, that is, $s_1=\lceil r / 2 \rceil$ and $s_2=\lceil (r+1)/ 2 \rceil$. Assuming $h=1/n$, $n\in\mathbb{N}$, we show the results in Tables~\ref{order_smooth_1} to~\ref{order_smooth_4}, where the error is $|P(0)-f_k(0)|$. Here $P(0)$~denotes the corresponding reconstruction from pointwise values to pointwise values at $x=0$. 

According to the tables, the numerical results are consistent with the theoretical considerations. Specifically,
 Table~\ref{order_smooth_1} shows that all  
  methods exhibit  the same order of accuracy for all  cases, dropping to second order when $k=1$, 
   in agreement with the previous theoretical considerations.  In Table~\ref{order_smooth_2} the gradual loss of accuracy as $k$~increases can be clearly seen  for WENO5 and CWENO5-LPR, whereas CWENO5 keeps the optimal accuracy order for $k<2$.  An analogous behavior can be seen in 
    Tables~\ref{order_smooth_3} and \ref{order_smooth_4}, corresponding to orders
    7 and 9, where WENO and CWENO-LPR gradually lose accuracy as $k$ increases, for $k=3$ and $4$, respectively, as expected.
 
\subsubsection*{Example 2: Discontinuous data analysis}

\begin{table}
 \addtolength{\tabcolsep}{-4.5pt}
  \centering
  \begin{tabular}{|c|c|c|c|c|c|c|c|c|c|c|c|c|}
    \hline
    & \multicolumn{2}{c|}{WENO3} & \multicolumn{2}{c|}{CWENO3-} & \multicolumn{2}{c|}{CWENO3} & \multicolumn{2}{c|}{WENO5} & \multicolumn{2}{c|}{CWENO5-} & \multicolumn{2}{c|}{CWENO5} \\
     & \multicolumn{2}{c|}{} & \multicolumn{2}{c|}{LPR} & \multicolumn{2}{c|}{} & \multicolumn{2}{c|}{} & \multicolumn{2}{c|}{LPR} & \multicolumn{2}{c|}{} \\
    \hline
    $n$ & Error & Ord. & Error & Ord. & Error & Ord. & Error & Ord. & Error & Ord. & Error & Ord. \\
    \hline
    40 & 2.68e-04 & --- & 4.31e-05 & --- & 7.65e-04 & --- & 4.58e-06 & --- & 4.73e-06 & --- & 4.47e-06 & --- \\
    \hline
    80 & 7.22e-05 & 1.89 & 9.16e-06 & 2.23 & 2.02e-04 & 1.92 & 5.91e-07 & 2.95 & 6.01e-07 & 2.98 & 5.84e-07 & 2.94 \\
    \hline
    160 & 1.87e-05 & 1.95 & 2.08e-06 & 2.14 & 5.21e-05 & 1.96 & 7.51e-08 & 2.98 & 7.57e-08 & 2.99 & 7.46e-08 & 2.97 \\
    \hline
    320 & 4.77e-06 & 1.97 & 4.93e-07 & 2.08 & 1.32e-05 & 1.98 & 9.46e-09 & 2.99 & 9.50e-09 & 2.99 & 9.43e-09 & 2.98 \\
    \hline
    640 & 1.20e-06 & 1.99 & 1.20e-07 & 2.04 & 3.32e-06 & 1.99 & 1.19e-09 & 2.99 & 1.19e-09 & 3.00 & 1.19e-09 & 2.99 \\
    \hline
    1280 & 3.03e-07 & 1.99 & 2.95e-08 & 2.02 & 8.34e-07 & 1.99 & 1.49e-10 & 3.00 & 1.49e-10 & 3.00 & 1.49e-10 & 3.00 \\
    \hline
    2560 & 7.58e-08 & 2.00 & 7.32e-09 & 2.01 & 2.09e-07 & 2.00 & 1.86e-11 & 3.00 & 1.86e-11 & 3.00 & 1.86e-11 & 3.00 \\
    \hline
    5120 & 1.90e-08 & 2.00 & 1.82e-09 & 2.01 & 5.23e-08 & 2.00 & 2.33e-12 & 3.00 & 2.33e-12 & 3.00 & 2.33e-12 & 3.00 \\
    \hline
    10240 & 4.75e-09 & 2.00 & 4.55e-10 & 2.00 & 1.31e-08 & 2.00 & 2.91e-13 & 3.00 & 2.91e-13 & 3.00 & 2.91e-13 & 3.00 \\
    \hline
  \end{tabular}
  \caption{Example~2: errors for $r=1$ and $r=2$.}
  \label{order_discontinuous_1}
\end{table}
\begin{table}
\addtolength{\tabcolsep}{-4.5pt}
  \centering
  \begin{tabular}{|c|c|c|c|c|c|c|c|c|c|c|c|c|}
    \hline
     & \multicolumn{2}{c|}{WENO7} & \multicolumn{2}{c|}{CWENO7-} & \multicolumn{2}{c|}{CWENO7} & \multicolumn{2}{c|}{WENO9} & \multicolumn{2}{c|}{CWENO9-} & \multicolumn{2}{c|}{CWENO9} \\
      & \multicolumn{2}{c|}{} & \multicolumn{2}{c|}{LPR} & \multicolumn{2}{c|}{} & \multicolumn{2}{c|}{} & \multicolumn{2}{c|}{LPR} & \multicolumn{2}{c|}{} \\
    \hline
    $n$ & Error & Ord. & Error & Ord. & Error & Ord. & Error & Ord. & Error & Ord. & Error & Ord. \\
    \hline
    40 & 2.44e-08 & --- & 9.68e-08 & --- & 1.40e-07 & --- & 2.28e-09 & --- & 2.28e-09 & --- & 2.28e-09 & --- \\
    \hline
    80 & 1.73e-09 & 3.82 & 6.16e-09 & 3.97 & 9.10e-09 & 3.94 & 7.31e-11 & 4.96 & 7.32e-11 & 4.96 & 7.31e-11 & 4.96 \\
    \hline
    160 & 1.15e-10 & 3.91 & 3.89e-10 & 3.99 & 5.81e-10 & 3.97 & 2.32e-12 & 4.98 & 2.32e-12 & 4.98 & 2.32e-12 & 4.98 \\
    \hline
    320 & 7.41e-12 & 3.96 & 2.44e-11 & 3.99 & 3.67e-11 & 3.98 & 7.29e-14 & 4.99 & 7.29e-14 & 4.99 & 7.28e-14 & 4.99 \\
    \hline
    640 & 4.70e-13 & 3.98 & 1.53e-12 & 4.00 & 2.31e-12 & 3.99 & 2.28e-15 & 5.00 & 2.28e-15 & 5.00 & 2.28e-15 & 5.00 \\
    \hline
    1280 & 2.96e-14 & 3.99 & 9.57e-14 & 4.00 & 1.44e-13 & 4.00 & 7.15e-17 & 5.00 & 7.15e-17 & 5.00 & 7.15e-17 & 5.00 \\
    \hline
    2560 & 1.86e-15 & 3.99 & 5.98e-15 & 4.00 & 9.04e-15 & 4.00 & 2.24e-18 & 5.00 & 2.24e-18 & 5.00 & 2.24e-18 & 5.00 \\
    \hline
    5120 & 1.16e-16 & 4.00 & 3.74e-16 & 4.00 & 5.65e-16 & 4.00 & 6.99e-20 & 5.00 & 6.99e-20 & 5.00 & 6.99e-20 & 5.00 \\
    \hline
    10240 & 7.28e-18 & 4.00 & 2.34e-17 & 4.00 & 3.54e-17 & 4.00 & 2.19e-21 & 5.00 & 2.19e-21 & 5.00 & 2.19e-21 & 5.00 \\
    \hline
  \end{tabular}
  \caption{Example~2: errors  for $r=3$ and $r=4$.}
  \label{order_discontinuous_2}
\end{table}

We now consider  the function $g:\mathbb{R}\to\mathbb{R}$ given by
\[g(x)=
\begin{cases}
  \mathrm{e}^x & \text{for $x\leq0$,}  \\
  \mathrm{e}^{x+1} &  \text{for $x>0$,} 
\end{cases}
\]
and perform the accuracy tests for $1\leq r\leq 4$ with the same setup as in Example~1. Results are shown in Tables~\ref{order_discontinuous_1} and~\ref{order_discontinuous_2}. It can be seen that the obtained accuracy order is consistent with our theoretical considerations. From the results it can be concluded that the order in presence of discontinuities is the optimal through the indicated choices for the parameters~$s_1$ and~$s_2$. This feature is thus shared by the classical WENO schemes with the suitable choice of the parameter~$s$.

\subsubsection*{Example 3: Non-aligned stencil with smooth data}

\begin{table}
  \centering
  \begin{tabular}{|c|c|c|c|c|c|c|c|c|}
    \hline
     & \multicolumn{4}{c|}{CWENO3-LPR} & \multicolumn{4}{c|}{CWENO3} \\
    \hline
    $r=1$ & \multicolumn{2}{c|}{$k=0$} & \multicolumn{2}{c|}{$k=1$} & \multicolumn{2}{c|}{$k=0$} & \multicolumn{2}{c|}{$k=1$} \\
    \hline
    $n$ & Error & Order & Error & Order & Error & Order & Error & Order \\
    \hline
    40 & 1.89e-05 & --- & 8.35e-05 & --- & 1.59e-06 & --- & 6.82e-05 & --- \\
    \hline
    80 & 2.00e-06 & 3.24 & 2.10e-05 & 1.99 & 2.61e-07 & 2.61 & 1.72e-05 & 1.99 \\
    \hline
    160 & 2.27e-07 & 3.14 & 5.28e-06 & 1.99 & 3.64e-08 & 2.84 & 4.31e-06 & 1.99 \\
    \hline
    320 & 2.70e-08 & 3.08 & 1.32e-06 & 2.00 & 4.78e-09 & 2.93 & 1.08e-06 & 2.00 \\
    \hline
    640 & 3.28e-09 & 3.04 & 3.31e-07 & 2.00 & 6.12e-10 & 2.97 & 2.70e-07 & 2.00 \\
    \hline
    1280 & 4.04e-10 & 3.02 & 8.27e-08 & 2.00 & 7.73e-11 & 2.98 & 6.75e-08 & 2.00 \\
    \hline
    2560 & 5.02e-11 & 3.01 & 2.07e-08 & 2.00 & 9.72e-12 & 2.99 & 1.69e-08 & 2.00 \\
    \hline
    5120 & 6.25e-12 & 3.01 & 5.17e-09 & 2.00 & 1.22e-12 & 3.00 & 4.22e-09 & 2.00 \\
    \hline
    10240 & 7.80e-13 & 3.00 & 1.29e-09 & 2.00 & 1.53e-13 & 3.00 & 1.05e-09 & 2.00 \\
    \hline
  \end{tabular}
  \caption{Example~3: errors  of schemes CWENO3-LPR and CWENO3.}
  \label{order_smooth_displ_1}
\end{table}
\begin{table}
  \centering
  \begin{tabular}{|c|c|c|c|c|c|c|c|}
    \hline
    & $r=2$ & \multicolumn{2}{c|}{$k=0$} & \multicolumn{2}{c|}{$k=1$} & \multicolumn{2}{c|}{$k=2$} \\
    \cline{2-8}
    & $n$ & Error & Order & Error & Order & Error & Order \\
    \hline
    \multirow{10}{*}{\rotatebox[origin=c]{90}{CWENO5-LPR}}
    & 40 & 4.72e-08 & --- & 2.36e-06 & --- & 3.67e-06 & --- \\
    \cline{2-8}
    & 80 & 1.45e-09 & 5.02 & 1.04e-07 & 4.50 & 4.56e-07 & 3.01 \\
    \cline{2-8}
    & 160 & 4.50e-11 & 5.01 & 5.11e-09 & 4.35 & 5.68e-08 & 3.01 \\
    \cline{2-8}
    & 320 & 1.40e-12 & 5.00 & 2.75e-10 & 4.21 & 7.08e-09 & 3.00 \\
    \cline{2-8}
    & 640 & 4.37e-14 & 5.00 & 1.59e-11 & 4.12 & 8.84e-10 & 3.00 \\
    \cline{2-8}
    & 1280 & 1.37e-15 & 5.00 & 9.49e-13 & 4.06 & 1.10e-10 & 3.00 \\
    \cline{2-8}
    & 2560 & 4.26e-17 & 5.00 & 5.80e-14 & 4.03 & 1.38e-11 & 3.00 \\
    \cline{2-8}
    & 5120 & 1.33e-18 & 5.00 & 3.58e-15 & 4.02 & 1.73e-12 & 3.00 \\
    \cline{2-8}
    & 10240 & 4.16e-20 & 5.00 & 2.23e-16 & 4.01 & 2.16e-13 & 3.00 \\
    \cline{2-8}
    & 20480 & 1.30e-21 & 5.00 & 1.39e-17 & 4.00 & 2.70e-14 & 3.00 \\
    \hline
    \multirow{10}{*}{\rotatebox[origin=c]{90}{CWENO5}}
    & 40 & 4.65e-10 & --- & 2.62e-07 & --- & 3.56e-06 & --- \\
    \cline{2-8}
    & 80 & 1.43e-11 & 5.02 & 8.03e-09 & 5.03 & 4.45e-07 & 3.00 \\
    \cline{2-8}
    & 160 & 4.47e-13 & 5.00 & 2.48e-10 & 5.02 & 5.56e-08 & 3.00 \\
    \cline{2-8}
    & 320 & 1.40e-14 & 5.00 & 7.68e-12 & 5.01 & 6.94e-09 & 3.00 \\
    \cline{2-8}
    & 640 & 4.37e-16 & 5.00 & 2.39e-13 & 5.01 & 8.68e-10 & 3.00 \\
    \cline{2-8}
    & 1280 & 1.37e-17 & 5.00 & 7.45e-15 & 5.00 & 1.08e-10 & 3.00 \\
    \cline{2-8}
    & 2560 & 4.27e-19 & 5.00 & 2.32e-16 & 5.00 & 1.36e-11 & 3.00 \\
    \cline{2-8}
    & 5120 & 1.34e-20 & 5.00 & 7.26e-18 & 5.00 & 1.69e-12 & 3.00 \\
    \cline{2-8}
    & 10240 & 4.17e-22 & 5.00 & 2.27e-19 & 5.00 & 2.12e-13 & 3.00 \\
    \cline{2-8}
    & 20480 & 1.30e-23 & 5.00 & 7.09e-21 & 5.00 & 2.65e-14 & 3.00 \\
    \hline
  \end{tabular}
  \caption{Example~3: errors of schemes CWENO5-LPR and CWENO5.}
  \label{order_smooth_displ_2}
\end{table}
\begin{table}
  \setlength\tabcolsep{4pt}
  \centering
  \begin{tabular}{|c|c|c|c|c|c|c|c|c|c|}
    \hline
    & $r=3$ & \multicolumn{2}{c|}{$k=0$} & \multicolumn{2}{c|}{$k=1$} & \multicolumn{2}{c|}{$k=2$} & \multicolumn{2}{c|}{$k=3$} \\
    \cline{2-10}
    & $n$ & Error & Order & Error & Order & Error & Order & Error & Order \\
    \hline
    \multirow{19}{*}{\rotatebox[origin=c]{90}{CWENO7-LPR}}
    & 40 & 1.59e-11 & --- & 1.13e-08 & --- & 1.77e-07 & --- & 1.46e-07 & --- \\
    \cline{2-10}
    & 80 & 1.08e-13 & 7.20 & 1.79e-10 & 5.98 & 1.65e-09 & 6.75 & 9.06e-09 & 4.01 \\
    \cline{2-10}
    & 160 & 7.80e-16 & 7.11 & 2.82e-12 & 5.99 & 7.23e-12 & 7.83 & 5.63e-10 & 4.01 \\
    \cline{2-10}
    & 320 & 5.84e-18 & 7.06 & 4.43e-14 & 5.99 & 1.07e-12 & 2.76 & 3.51e-11 & 4.01 \\
    \cline{2-10}
    & 640 & 4.46e-20 & 7.03 & 6.94e-16 & 6.00 & 4.57e-14 & 4.55 & 2.19e-12 & 4.00 \\
    \cline{2-10}
    & 1280 & 3.45e-22 & 7.02 & 1.08e-17 & 6.00 & 1.61e-15 & 4.82 & 1.37e-13 & 4.00 \\
    \cline{2-10}
    & 2560 & 2.68e-24 & 7.01 & 1.70e-19 & 6.00 & 5.33e-17 & 4.92 & 8.53e-15 & 4.00 \\
    \cline{2-10}
    & 5120 & 2.09e-26 & 7.00 & 2.65e-21 & 6.00 & 1.71e-18 & 4.96 & 5.33e-16 & 4.00 \\
    \cline{2-10}
    & 10240 & 1.63e-28 & 7.00 & 4.14e-23 & 6.00 & 5.41e-20 & 4.98 & 3.33e-17 & 4.00 \\
    \cline{2-10}
    & 20480 & 1.27e-30 & 7.00 & 6.47e-25 & 6.00 & 1.70e-21 & 4.99 & 2.08e-18 & 4.00 \\
    \cline{2-10}
    & 40960 & 9.93e-33 & 7.00 & 1.01e-26 & 6.00 & 5.33e-23 & 5.00 & 1.30e-19 & 4.00 \\
    \cline{2-10}
    & 81920 & 7.76e-35 & 7.00 & 1.58e-28 & 6.00 & 1.67e-24 & 5.00 & 8.13e-21 & 4.00 \\
    \cline{2-10}
    & 163840 & 6.06e-37 & 7.00 & 2.47e-30 & 6.00 & 5.22e-26 & 5.00 & 5.08e-22 & 4.00 \\
    \cline{2-10}
    & 327680 & 4.73e-39 & 7.00 & 3.86e-32 & 6.00 & 1.63e-27 & 5.00 & 3.18e-23 & 4.00 \\
    \cline{2-10}
    & 655360 & 3.70e-41 & 7.00 & 6.03e-34 & 6.00 & 5.10e-29 & 5.00 & 1.99e-24 & 4.00 \\
    \cline{2-10}
    & 1310720 & 2.89e-43 & 7.00 & 9.42e-36 & 6.00 & 1.59e-30 & 5.00 & 1.24e-25 & 4.00 \\
    \cline{2-10}
    & 2621440 & 2.26e-45 & 7.00 & 1.47e-37 & 6.00 & 4.98e-32 & 5.00 & 7.76e-27 & 4.00 \\
    \cline{2-10}
    & 5242880 & 1.76e-47 & 7.00 & 2.30e-39 & 6.00 & 1.56e-33 & 5.00 & 4.85e-28 & 4.00 \\
    \cline{2-10}
    & 10485760 & 1.38e-49 & 7.00 & 3.59e-41 & 6.00 & 4.86e-35 & 5.00 & 3.03e-29 & 4.00 \\
    \hline
    \multirow{19}{*}{\rotatebox[origin=c]{90}{CWENO7}}
    & 40 & 7.92e-14 & --- & 4.74e-13 & --- & 1.14e-09 & --- & 1.00e-07 & --- \\
    \cline{2-10}
    & 80 & 6.24e-16 & 6.99 & 3.74e-15 & 6.98 & 4.67e-12 & 7.94 & 6.48e-09 & 3.95 \\
    \cline{2-10}
    & 160 & 4.90e-18 & 6.99 & 2.94e-17 & 6.99 & 1.83e-14 & 8.00 & 4.11e-10 & 3.98 \\
    \cline{2-10}
    & 320 & 3.84e-20 & 7.00 & 2.30e-19 & 7.00 & 7.07e-17 & 8.01 & 2.59e-11 & 3.99 \\
    \cline{2-10}
    & 640 & 3.00e-22 & 7.00 & 1.80e-21 & 7.00 & 2.71e-19 & 8.03 & 1.62e-12 & 3.99 \\
    \cline{2-10}
    & 1280 & 2.35e-24 & 7.00 & 1.41e-23 & 7.00 & 1.02e-21 & 8.05 & 1.02e-13 & 4.00 \\
    \cline{2-10}
    & 2560 & 1.83e-26 & 7.00 & 1.10e-25 & 7.00 & 3.71e-24 & 8.10 & 6.36e-15 & 4.00 \\
    \cline{2-10}
    & 5120 & 1.43e-28 & 7.00 & 8.60e-28 & 7.00 & 1.23e-26 & 8.23 & 3.98e-16 & 4.00 \\
    \cline{2-10}
    & 10240 & 1.12e-30 & 7.00 & 6.72e-30 & 7.00 & 3.14e-29 & 8.62 & 2.49e-17 & 4.00 \\
    \cline{2-10}
    & 20480 & 8.75e-33 & 7.00 & 5.25e-32 & 7.00 & 8.65e-33 & 11.82 & 1.55e-18 & 4.00 \\
    \cline{2-10}
    & 40960 & 6.83e-35 & 7.00 & 4.10e-34 & 7.00 & 1.06e-33 & 3.03 & 9.72e-20 & 4.00 \\
    \cline{2-10}
    & 81920 & 5.34e-37 & 7.00 & 3.20e-36 & 7.00 & 1.21e-35 & 6.45 & 6.07e-21 & 4.00 \\
    \cline{2-10}
    & 163840 & 4.17e-39 & 7.00 & 2.50e-38 & 7.00 & 1.10e-37 & 6.79 & 3.80e-22 & 4.00 \\
    \cline{2-10}
    & 327680 & 3.26e-41 & 7.00 & 1.96e-40 & 7.00 & 9.18e-40 & 6.90 & 2.37e-23 & 4.00 \\
    \cline{2-10}
    & 655360 & 2.55e-43 & 7.00 & 1.53e-42 & 7.00 & 7.41e-42 & 6.95 & 1.48e-24 & 4.00 \\
    \cline{2-10}
    & 1310720 & 1.99e-45 & 7.00 & 1.19e-44 & 7.00 & 5.88e-44 & 6.98 & 9.27e-26 & 4.00 \\
    \cline{2-10}
    & 2621440 & 1.55e-47 & 7.00 & 9.32e-47 & 7.00 & 4.63e-46 & 6.99 & 5.79e-27 & 4.00 \\
    \cline{2-10}
    & 5242880 & 1.21e-49 & 7.00 & 7.28e-49 & 7.00 & 3.63e-48 & 6.99 & 3.62e-28 & 4.00 \\
    \cline{2-10}
    & 10485760 & 9.48e-52 & 7.00 & 5.69e-51 & 7.00 & 2.84e-50 & 7.00 & 2.26e-29 & 4.00 \\
    \hline
  \end{tabular}
  \caption{Example~3: errors of schemes CWENO7-LPR and CWENO7.}
  \label{order_smooth_displ_3}
\end{table}
\begin{table}
 \addtolength{\tabcolsep}{-3.5pt} 
  \centering
  \begin{tabular}{|c|c|c|c|c|c|c|c|c|c|c|c|}
    \hline
    & $r=4$ & \multicolumn{2}{c|}{$k=0$} & \multicolumn{2}{c|}{$k=1$} & \multicolumn{2}{c|}{$k=2$} & \multicolumn{2}{c|}{$k=3$} & \multicolumn{2}{c|}{$k=4$} \\
    \cline{2-12}
    & $n$ & Error & Ord. & Error & Ord. & Error & Ord. & Error & Ord. & Error & Ord. \\
    \hline
    \multirow{9}{*}{\rotatebox[origin=c]{90}{CWENO9-LPR}}
    & 40 & 2.79e-14 & --- & 7.07e-12 & --- & 5.44e-09 & --- & 6.98e-09 & --- & 9.25e-09 & --- \\
    \cline{2-12}
    & 80 & 5.38e-17 & 9.02 & 2.32e-14 & 8.25 & 4.48e-11 & 6.92 & 9.15e-10 & 2.93 & 2.85e-10 & 5.02 \\
    \cline{2-12}
    & 160 & 1.05e-19 & 9.01 & 8.16e-17 & 8.15 & 3.57e-13 & 6.97 & 2.05e-11 & 5.48 & 8.84e-12 & 5.01 \\
    \cline{2-12}
    & 320 & 2.04e-22 & 9.00 & 3.01e-19 & 8.08 & 2.82e-15 & 6.99 & 3.48e-13 & 5.88 & 2.75e-13 & 5.01 \\
    \cline{2-12}
    & 640 & 3.98e-25 & 9.00 & 1.14e-21 & 8.04 & 2.21e-17 & 6.99 & 5.56e-15 & 5.97 & 8.59e-15 & 5.00 \\
    \cline{2-12}
    & 1280 & 7.78e-28 & 9.00 & 4.39e-24 & 8.02 & 1.73e-19 & 7.00 & 8.74e-17 & 5.99 & 2.68e-16 & 5.00 \\
    \cline{2-12}
    & 2560 & 1.52e-30 & 9.00 & 1.70e-26 & 8.01 & 1.35e-21 & 7.00 & 1.37e-18 & 6.00 & 8.37e-18 & 5.00 \\
    \cline{2-12}
    & 5120 & 2.97e-33 & 9.00 & 6.62e-29 & 8.01 & 1.06e-23 & 7.00 & 2.14e-20 & 6.00 & 2.62e-19 & 5.00 \\
    \cline{2-12}
    & 10240 & 5.79e-36 & 9.00 & 2.58e-31 & 8.00 & 8.26e-26 & 7.00 & 3.34e-22 & 6.00 & 8.18e-21 & 5.00 \\
    \hline
    \multirow{9}{*}{\rotatebox[origin=c]{90}{CWENO9}}
    & 40 & 1.36e-17 & --- & 1.09e-16 & --- & 1.25e-15 & --- & 1.18e-09 & --- & 9.07e-09 & --- \\
    \cline{2-12}
    & 80 & 2.69e-20 & 8.99 & 2.15e-19 & 8.99 & 1.56e-18 & 9.64 & 2.85e-12 & 8.69 & 2.81e-10 & 5.01 \\
    \cline{2-12}
    & 160 & 5.28e-23 & 8.99 & 4.22e-22 & 8.99 & 2.96e-21 & 9.05 & 5.82e-15 & 8.94 & 8.76e-12 & 5.01 \\
    \cline{2-12}
    & 320 & 1.03e-25 & 9.00 & 8.26e-25 & 9.00 & 5.78e-24 & 9.00 & 1.14e-17 & 9.00 & 2.73e-13 & 5.00 \\
    \cline{2-12}
    & 640 & 2.02e-28 & 9.00 & 1.62e-27 & 9.00 & 1.13e-26 & 9.00 & 2.22e-20 & 9.01 & 8.52e-15 & 5.00 \\
    \cline{2-12}
    & 1280 & 3.95e-31 & 9.00 & 3.16e-30 & 9.00 & 2.21e-29 & 9.00 & 4.32e-23 & 9.00 & 2.66e-16 & 5.00 \\
    \cline{2-12}
    & 2560 & 7.71e-34 & 9.00 & 6.17e-33 & 9.00 & 4.32e-32 & 9.00 & 8.42e-26 & 9.00 & 8.31e-18 & 5.00 \\
    \cline{2-12}
    & 5120 & 1.51e-36 & 9.00 & 1.21e-35 & 9.00 & 8.44e-35 & 9.00 & 1.64e-28 & 9.00 & 2.60e-19 & 5.00 \\
    \cline{2-12}
    & 10240 & 2.94e-39 & 9.00 & 2.35e-38 & 9.00 & 1.65e-37 & 9.00 & 3.21e-31 & 9.00 & 8.11e-21 & 5.00 \\
    \hline
  \end{tabular}
  \caption{Example~3: errors of schemes CWENO9-LPR and CWENO9.}
  \label{order_smooth_displ_4}
\end{table}

We now consider again the setup as in Example 1, with the difference that now the stencil is based on the 
 grid points $x_j= (j- 3/4 )h$, $-r\leq j\leq r$,  which are based on choosing the displacement parameter $\tau=3/4$. 
Since we want to interpolate at $x=0$, in this case the stencil is displaced with respect to $x$, and thus the classical WENO procedure cannot provide a satisfactory procedure to solve this problem, whereas both CWENO-LPR and our CWENO approach are able to by just choosing, for instance, the subweights based on the uniform ideal weights, as it can be seen in Tables~\ref{order_smooth_displ_1} to~\ref{order_smooth_displ_4}.
The results of each of these tables are consistent with our theoretical considerations, since
Table~\ref{order_smooth_displ_1}, \ref{order_smooth_displ_2}, \ref{order_smooth_displ_3}, and 
  \ref{order_smooth_displ_4}, respectively, indicate that the optimal order is attained for 
 $k<1$, $k<2$, $k<3$, and $k<4$,  respectively. 
The results show that our scheme is also suitable for problems where
the reconstruction point is not centered and attains the optimal
order as in the centered case.

\subsubsection*{Example 4: Non-aligned stencil with discontinuous data}

\begin{table}
  \setlength\tabcolsep{5pt}
  \centering
  \begin{tabular}{|c|c|c|c|c|c|c|c|c|c|}
    \hline
    & & \multicolumn{2}{c|}{$r=1$} & \multicolumn{2}{c|}{$r=2$} & \multicolumn{2}{c|}{$r=3$} & \multicolumn{2}{c|}{$r=4$} \\
    \cline{2-10}
    & $n$ & Error & Order & Error & Order & Error & Order & Error & Order \\
    \hline
    \multirow{8}{*}{\rotatebox[origin=c]{90}{CWENO-LPR}}
    & 40 & 7.12e-04 & --- & 8.88e-06 & --- & 1.34e-08 & --- & 4.94e-09 & --- \\
    \cline{2-10}
    & 80 & 1.91e-04 & 1.90 & 1.14e-06 & 2.96 & 8.07e-10 & 4.05 & 1.59e-10 & 4.96 \\
    \cline{2-10}
    & 160 & 4.94e-05 & 1.95 & 1.45e-07 & 2.98 & 4.94e-11 & 4.03 & 5.04e-12 & 4.98 \\
    \cline{2-10}
    & 320 & 1.26e-05 & 1.98 & 1.82e-08 & 2.99 & 3.05e-12 & 4.02 & 1.59e-13 & 4.99 \\
    \cline{2-10}
    & 640 & 3.17e-06 & 1.99 & 2.29e-09 & 2.99 & 1.89e-13 & 4.01 & 4.97e-15 & 4.99 \\
    \cline{2-10}
    & 1280 & 7.95e-07 & 1.99 & 2.86e-10 & 3.00 & 1.18e-14 & 4.00 & 1.56e-16 & 5.00 \\
    \cline{2-10}
    & 2560 & 1.99e-07 & 2.00 & 3.58e-11 & 3.00 & 7.36e-16 & 4.00 & 4.87e-18 & 5.00 \\
    \cline{2-10}
    & 5120 & 4.98e-08 & 2.00 & 4.48e-12 & 3.00 & 4.60e-17 & 4.00 & 1.52e-19 & 5.00 \\
    \hline
    \multirow{8}{*}{\rotatebox[origin=c]{90}{CWENO}}
    & 40 & 1.13e-04 & --- & 9.07e-06 & --- & 2.00e-07 & --- & 4.94e-09 & --- \\
    \cline{2-10}
    & 80 & 2.58e-05 & 2.13 & 1.15e-06 & 2.97 & 1.28e-08 & 3.97 & 1.59e-10 & 4.96 \\
    \cline{2-10}
    & 160 & 6.14e-06 & 2.07 & 1.46e-07 & 2.99 & 8.07e-10 & 3.98 & 5.04e-12 & 4.98 \\
    \cline{2-10}
    & 320 & 1.49e-06 & 2.04 & 1.83e-08 & 2.99 & 5.07e-11 & 3.99 & 1.59e-13 & 4.99 \\
    \cline{2-10}
    & 640 & 3.69e-07 & 2.02 & 2.29e-09 & 3.00 & 3.18e-12 & 4.00 & 4.97e-15 & 4.99 \\
    \cline{2-10}
    & 1280 & 9.15e-08 & 2.01 & 2.87e-10 & 3.00 & 1.99e-13 & 4.00 & 1.56e-16 & 5.00 \\
    \cline{2-10}
    & 2560 & 2.28e-08 & 2.01 & 3.58e-11 & 3.00 & 1.24e-14 & 4.00 & 4.87e-18 & 5.00 \\
    \cline{2-10}
    & 5120 & 5.69e-09 & 2.00 & 4.48e-12 & 3.00 & 7.78e-16 & 4.00 & 1.52e-19 & 5.00 \\
    \hline
  \end{tabular}
  \caption{Example~4: errors of  
   schemes~CWENO($2r+1$)-LPR and~CWENO($2r+1$).} 
  \label{order_discontinuous_displ}
\end{table}

Now, we change our setup to the one defined in Example 2, with the stencil arrangement and subweights of Example 3. Results are shown in Table \ref{order_discontinuous_displ}. It can be concluded that in the case of non-centered reconstruction points the order is also the optimal in presence of discontinuities in the data.

\subsection{Experiments involving the numerical solution of conservation laws}\label{clex}

\subsubsection*{Example 5. 1D Euler equations: Shu-Osher problem}

\begin{figure}[t]
  \centering
  \begin{tabular}{cc}
    \includegraphics[width=0.47\textwidth]{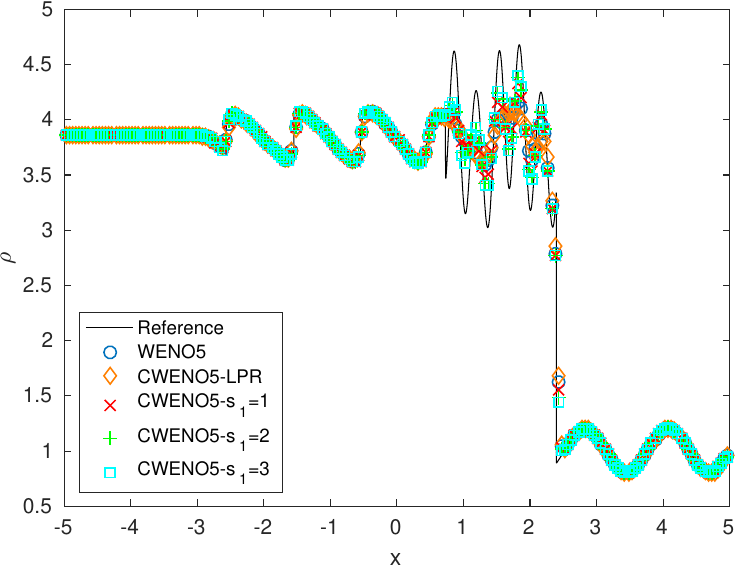} & \includegraphics[width=0.47\textwidth, height=0.22\textheight]{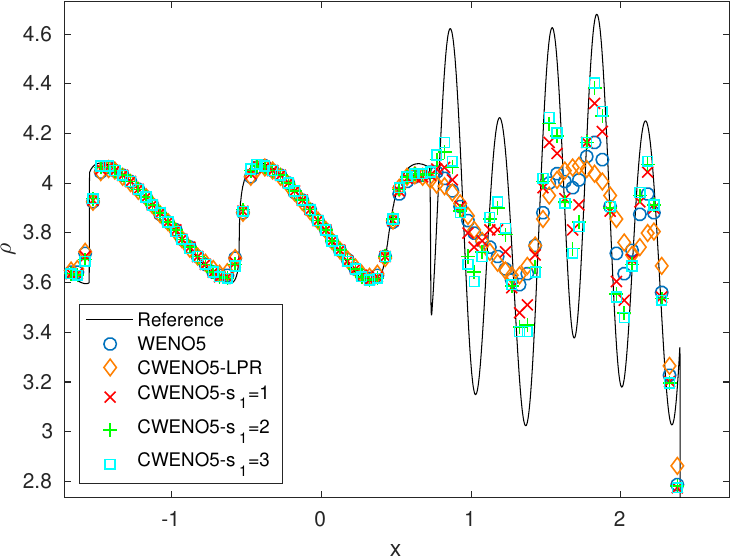} \\
    (a) Numerical solution overview & (b) Enlarged view (I) \\[3ex]
    \includegraphics[width=0.47\textwidth]{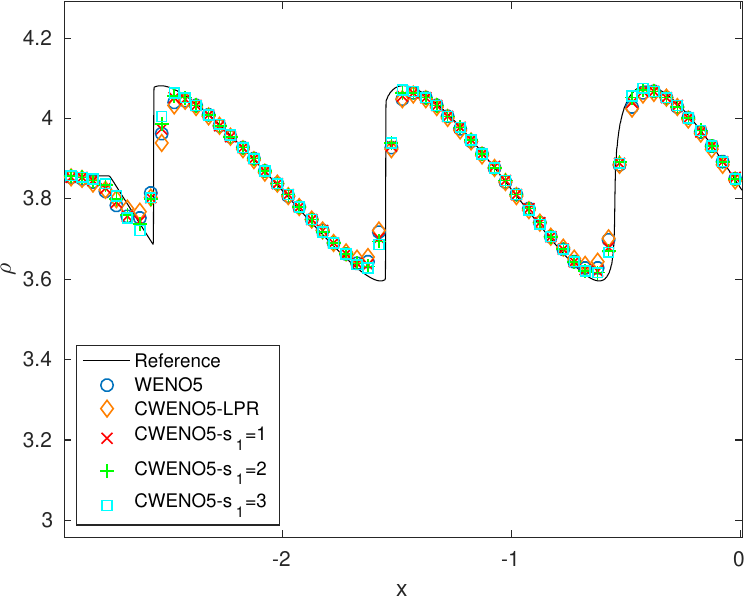} & \includegraphics[width=0.47\textwidth]{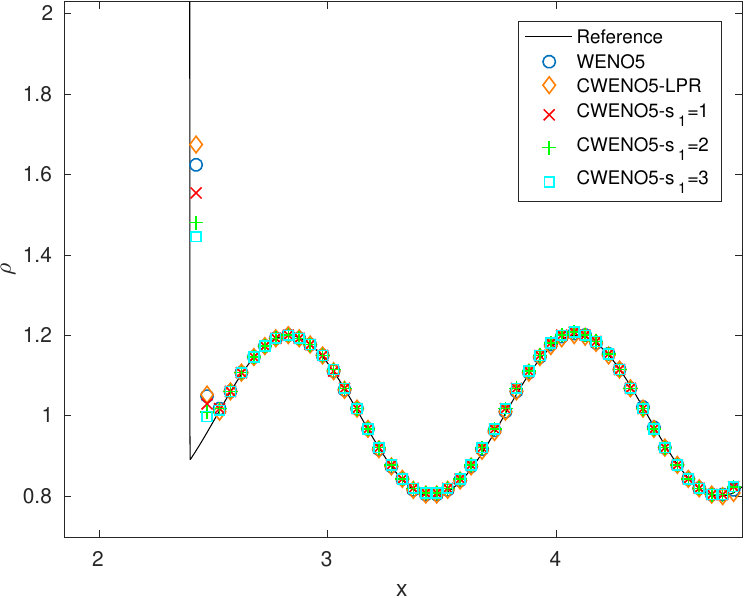} \\
    (c) Enlarged view (II) & (d) Enlarged view (III)
  \end{tabular}
  \caption{Example 5: density field, discretization with $200$ points, $s_2=2$, $T=1.8$.}
  \label{shuosher}
\end{figure}

\begin{figure}[t]
  \centering
  \includegraphics[width=0.55\textwidth]{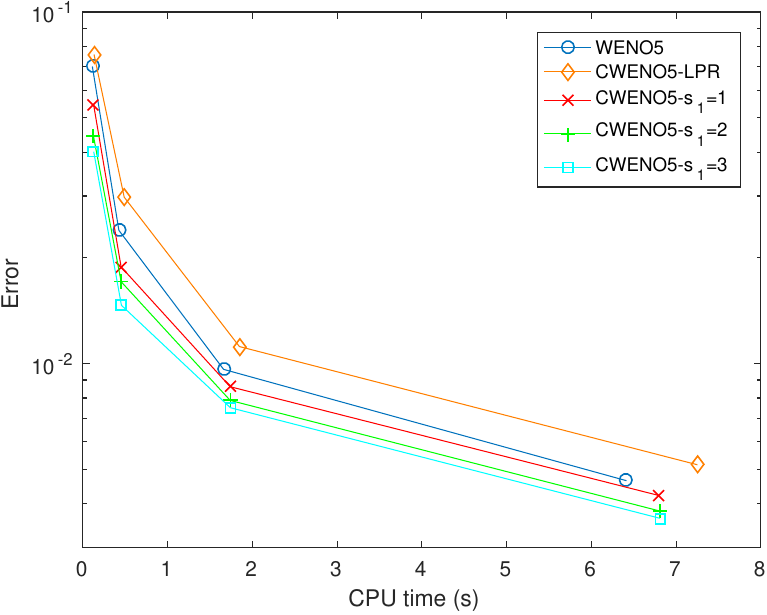}
  \caption{Example 5: efficiency plot.}
  \label{shuosher_cpu}
\end{figure}

The 1D Euler equations for gas dynamics are given by \eqref{hcl} for $\nu=3$ and $d=1$ with 
 $\boldsymbol{u} = ( \rho, \rho v, E)^{\mathrm{T}}$ and $\boldsymbol{f} (\boldsymbol{u}) = 
  \boldsymbol{f}^1  (\boldsymbol{u}) = (\rho v, 
        p+\rho v^2, 
        v(E+p))^{\mathrm{T}}$,  where $\rho$ is the density, $v$ is the velocity and  $E$ is the specific energy of the system. The variable~$p$ stands for the pressure and is given by the equation of state
$$p=\left(\gamma-1\right)\left(E-\frac{1}{2}\rho v^2\right),$$
where $\gamma$ is the adiabatic constant that will be taken as $\gamma =1.4$. 
We now consider the interaction with a Mach 3 shock and a sine wave. The spatial domain is now given by $\Omega:=(-5,5)\ni x_1 =: x$, with the initial condition 
\begin{align*} 
(\rho,v,p) (x, 0) = 
\begin{cases}
  (3.857143, 
    2.629369, 
    10.33333) & \text{if  $x\leq-4$,}  \\
   (1.0+0.2\sin(5x), 0, 1) 
    & \text{if $x>-4$,} 
\end{cases} \end{align*} 
with left inflow and right outflow boundary conditions. 

We run a simulation until $T=1.8$ and compare the results obtained with the classical WENO5 scheme, the CWENO-LPR method and the proposed CWENO schemes for $1\leq s_1\leq 3$ (let us recall that the smallest parameter to achieve fifth-order accuracy in this case is $r_1=1$), $s_2=2$ (the smallest parameter to achieve the optimal accuracy in presence of discontinuities), using the subweights based on the classical WENO ideal weights, $n=200$ cells, $\textnormal{CFL}=0.5$,  and a reference solution computed with $16000$ grid points. The results are shown in Figure~\ref{shuosher} for the density field.
 One can conclude that the new CWENO schemes capture better both the smooth extrema and the discontinuities in the numerical solution, yielding better results as the parameter $s_1$ increases (namely, when the global average weight involving the spatial reconstructions is closest to $1$, thus increasing the impact of the full degree reconstruction polynomial on the reconstruction).
Finally, in order to stress the performance of our new scheme, we also show a comparison between the schemes involving the CPU time versus the error in $L^1$-norm, which is shown in Figure~\ref{shuosher_cpu}.

\subsubsection*{Example 6. 2D Euler equations: double Mach reflection}

The two-dimensional Euler equations for inviscid gas dynamics
are  given by \eqref{hcl} for $\nu=4$ and $d=2$, where  for  $x=x_1$ and  $y=x_2$,  
 we set  \begin{align*} 
& \boldsymbol{u} =  \begin{pmatrix} 
    \rho \\
    \rho v^x \\
    \rho v^y \\
    E \end{pmatrix}, \quad 
  \boldsymbol{f}^1 (\boldsymbol{u})= \begin{pmatrix}  
    \rho v^x \\
    p+\rho (v^x)^2 \\
    \rho v^xv^y \\
    v^x(E+p) 
  \end{pmatrix},  \quad   \boldsymbol{f}^2 (\boldsymbol{u})= \begin{pmatrix} 
    \rho v^y \\
    \rho v^xv^y \\
    p+\rho (v^y)^2 \\
    v^y(E+p) 
  \end{pmatrix}. 
\end{align*}
Here  $\rho$ is the density, $(v^x, v^y)$  is the velocity, $E$ is the specific energy, and  $p$  is the pressure  that  is given by the equation of state $$p=(\gamma-1)\left(E-\frac{1}{2}\rho((v^x)^2+(v^y)^2)\right),$$ where
 the adiabatic constant is again chosen as  $\gamma =1.4$.

This experiment uses these equations to model a vertical right-going Mach 10 shock colliding with an equilateral triangle. By symmetry, this is equivalent to a collision with a ramp with a slope of $30^{\circ}$ with respect to the horizontal line.

\begin{figure}[t]
  \centering
  \setlength{\tabcolsep}{1cm}
  \begin{tabular}{cc}
    \includegraphics[width=0.35\textwidth]{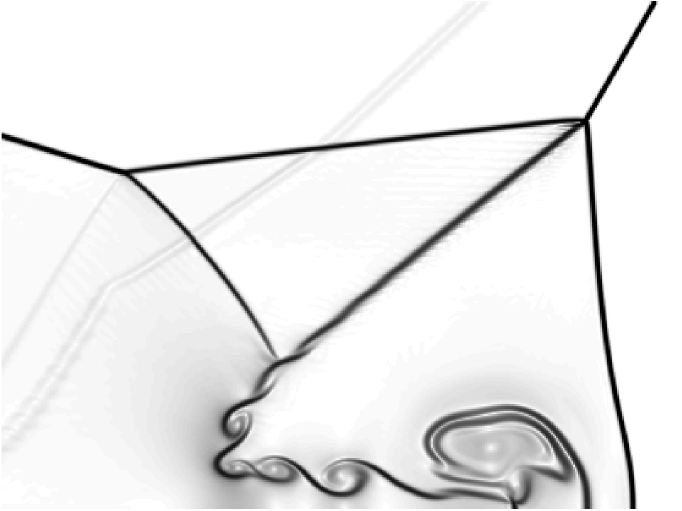} & \includegraphics[width=0.35\textwidth]{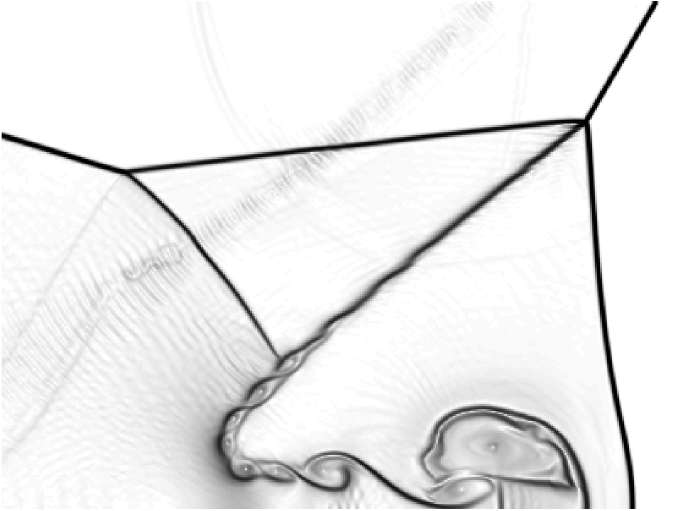} \\
    (a) WENO5, $2048\times512$ & (b) CWENO5-LPR, $2048\times512$ \\
    \includegraphics[width=0.35\textwidth]{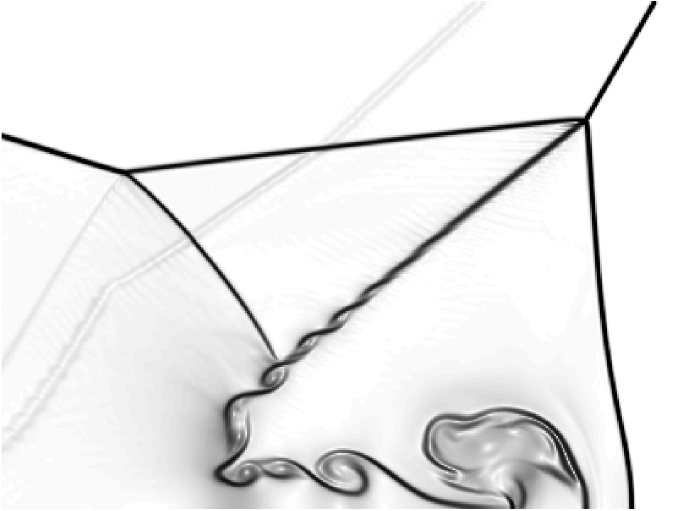} & \includegraphics[width=0.35\textwidth]{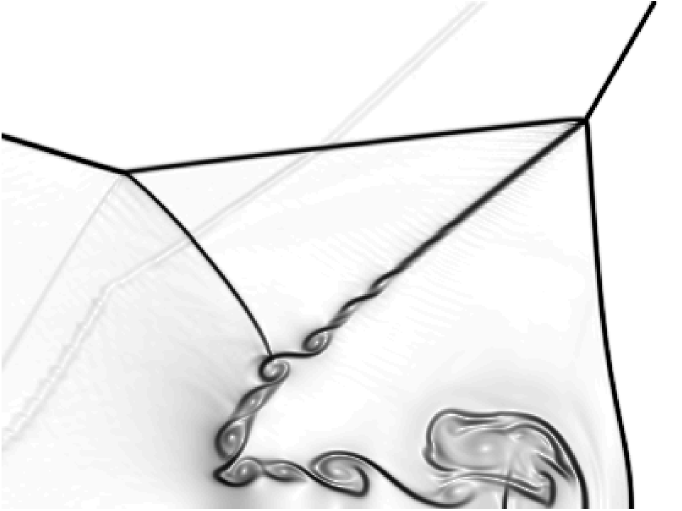} \\
    (c) CWENO5, $s_1=1$, $2048\times512$ & (d) CWENO5, $s_1=2$, $2048\times512$ \\
    \includegraphics[width=0.35\textwidth]{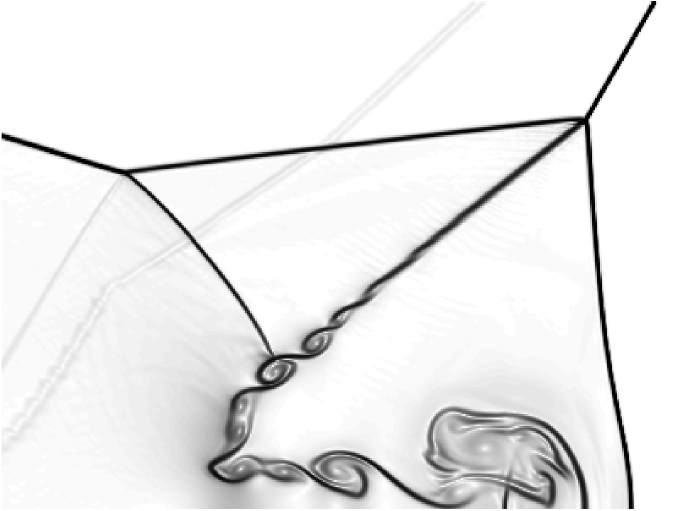} & \includegraphics[width=0.35\textwidth]{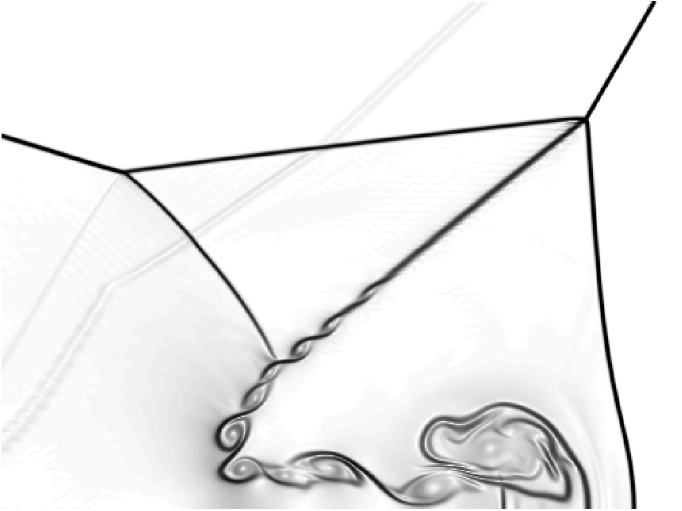}  \\
    (e) CWENO5, $s_1=3$, $2048\times512$ & (f) WENO5, $2560\times640$ \\
    \includegraphics[width=0.35\textwidth]{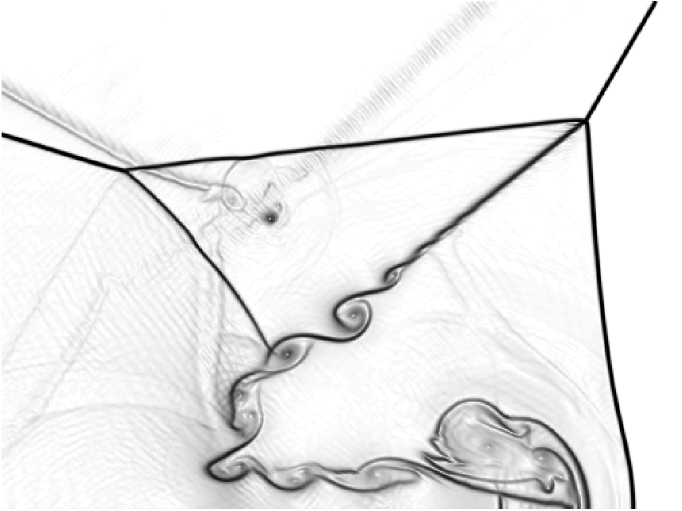} & \includegraphics[width=0.35\textwidth]{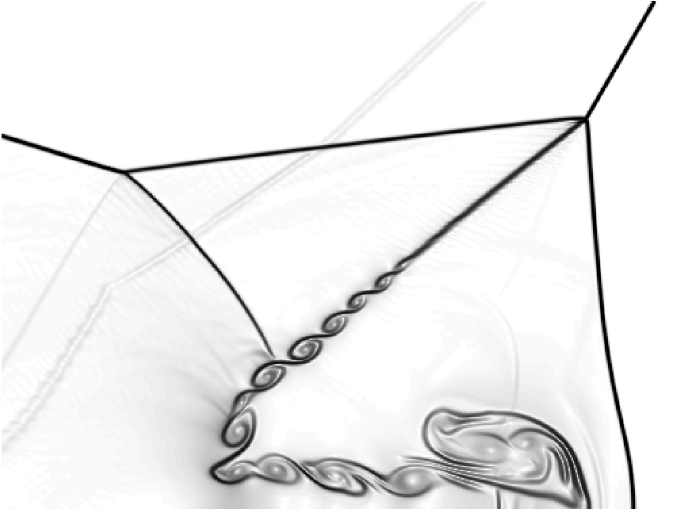}  \\
    (g) CWENO5-LPR, $2560\times640$ & (h) CWENO5, $s_1=1$, $2560\times640$ \\
    \includegraphics[width=0.35\textwidth]{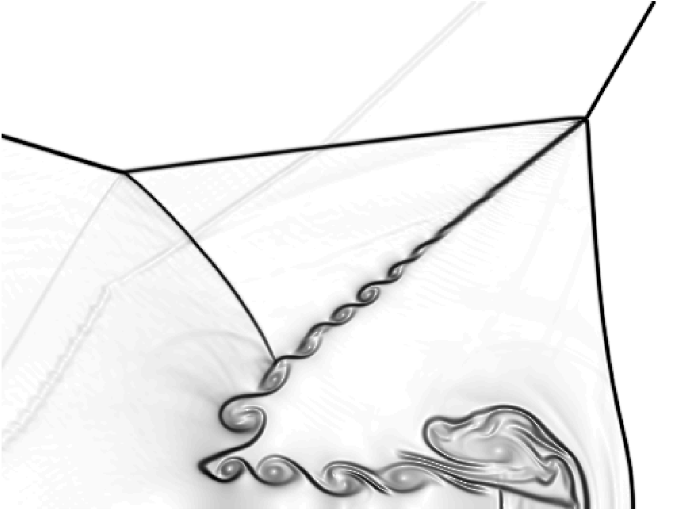} & \includegraphics[width=0.35\textwidth]{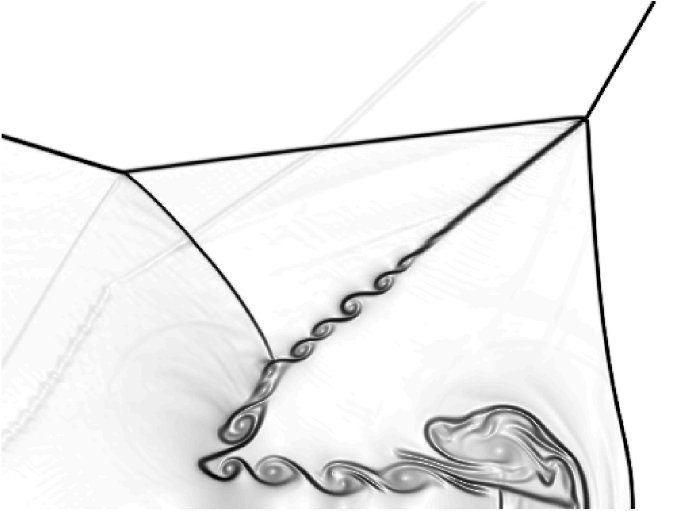}  \\
    (i) CWENO5, $s_1=2$, $2560\times640$ & (j) CWENO5, $s_1=3$, $2560\times640$
  \end{tabular}
  \caption{Example~6: enlarged view of turbulent zone (Schlieren plot).}
  \label{dmr}
\end{figure}

For  sake of simplicity, we consider the equivalent problem in a rectangle, consisting in a rotated shock, whose vertical angle is 
$30^{\circ}$. The domain is the rectangle $\Omega=[0,4]\times[0,1]$, whose initial conditions are
\begin{gather*}  (\rho,v^x,v^y,E) (x,y,0)=\begin{cases}
  \boldsymbol{c}_1= (\rho_1,v_1^x,v_1^y,E_1)    & \text{if $y\leq 1/4 +\tan(\pi/6)x$,} \\
  \boldsymbol{c}_2=  (\rho_2,v_2^x,v_2^y,E_2)     & \text{if $y > 1/4 +\tan(\pi/6)x$,}  
\end{cases} \\
    \boldsymbol{c}_1 =
    \bigl(8,8.25\cos(\pi/6),-8.25\sin(\pi/6),563.5\bigr), \quad 
    \boldsymbol{c}_2=    (1.4,0,0,2.5).
  \end{gather*} 
We impose inflow boundary conditions, with value $\boldsymbol{c}_1$, at the left side, $\{0\}\times[0,1]$, outflow boundary conditions both at $[0,1/4]\times\{0\}$ and $\{4\}\times[0,1]$, reflecting boundary conditions at  $(1/4,4]\times\{0\}$ and inflow boundary conditions at the upper side, $[0,4]\times\{1\}$, which mimics the shock at its actual traveling speed:
\begin{align*} 
 (\rho,v^x,v^y,E) (x,1,t)=\begin{cases}
  \boldsymbol{c}_1 & \text{if $x\leq 1/4 + (1+20t)/\sqrt{3}$,}  \\
  \boldsymbol{c}_2 & \text{if $x>1/4 + (1+20t)/\sqrt{3}$.}  
\end{cases} \end{align*} 
We run different simulations until $T=0.2$ both at a resolution of $2048\times512$ points and a resolution of $2560\times640$ points, shown in Figure \ref{dmr}, in both cases with $\textnormal{CFL}=0.4$ and involving the classical WENO5 scheme, the CWENO5-LPR method and our CWENO schemes with $1\leq s_1\leq 3$, $s_2=2$, and using the subweights based on the classical WENO ideal weights in the latter case.

From the results, it can be seen that our scheme captures better some features of a complex weak solution, such as turbulence and vorticity, than the corresponding WENO and CWENO-LPR counterparts. Moreover, and in qualitative terms, the results obtained with a CWENO5 scheme in a resolution of $2048\times512$ points are similar to those obtained through the classical WENO5 and CWENO schemes in a resolution of $2560\times640$ points. In addition, and as it can be observed in the figures, CWENO5-LPR seems to generate moderate spurious oscillations along the whole solution in a resolution of $2048\times512$ points. These oscillations are more pronounced in a resolution of $2560\times640$ points, which does not occur in the case of WENO and CWENO schemes. 
Finally, Table \ref{cpu} shows a comparison  of the computational cost for the corresponding $(2r+1)$-th order schemes, with $1\leq r\leq4$, for a resolution of $256\times64$ points.

\begin{table}
  \setlength\tabcolsep{4.5pt}
  \centering
  \begin{tabular}{|c|c|c|c|c|}
    \hline
    Schemes/Cost & $r=1$ & $r=2$ & $r=3$ & $r=4$ \\
    \hline
    WENO & 36.073758 & 56.399993 & 79.428691 & 106.035864 \\
    \hline
    CWENO-LPR & 41.929048 & 70.513481 & 98.687388 & 144.724150 \\
    \hline
    CWENO & 39.251961 & 63.518240 & 86.209761 & 110.557309 \\
    \hline
    Ratio CWENO-LPR/WENO & 1.1623 & 1.2502 & 1.2425 & 1.3649 \\
    \hline
    Ratio CWENO/WENO & 1.0881 & 1.1262 & 1.0854 & 1.0426 \\
    \hline
    Ratio CWENO/CWENO-LPR & 0.9362 & 0.9008 & 0.8736 & 0.7639 \\
    \hline
  \end{tabular}
  \caption{Example 6: efficiency table for a grid of  $256\times64$ points (cost in seconds).}
  \label{cpu}
\end{table}

\subsubsection*{Example 7. 2D Euler equations: Riemann problem}

We now consider a last experiment consisting in a Riemann problem for the 2D Euler equations on the domain $(0,1)\times(0,1)$. Riemann problems for 2D Euler equations were first studied in \cite{SchulzRinne}. The initial data is taken from \cite[Sect.~3, Config.~3]{KurganovTadmor}:
\begin{align*}
  \boldsymbol{u}(x, y, 0)=(\rho(x, y, 0), \rho(x, y, 0)v^x(x, y, 0), \rho(x, y, 0)v^y(x, y, 0), E(x, y, 0))
  \end{align*}
and
\[
\begin{pmatrix}
  \rho(x, y, 0)\\
  v^x(x, y, 0)\\
  v^y(x, y, 0)\\
  p(x, y, 0)
\end{pmatrix}^{\mathrm{T}} 
  =
\begin{cases}
  (1.5, 0, 0, 1.5) & \text{for $x>0.5$, $y>0.5$,}  \\
  (0.5323, 1.206, 0, 0.3) & \text{for  $x\leq0.5$, $y>0.5$,}  \\
  (0.138, 1.206, 1.206, 0.029) &  \text{for $x\leq0.5$, $y\leq0.5$,}  \\
  (0.5323, 0, 1.206, 0.3) &  \text{for $x>0.5$, $y\leq0.5$,} 
\end{cases}
\]
with the same equation of state as in the previous test.

\begin{figure}[t]
  \centering
  \setlength{\tabcolsep}{1cm}
  \begin{tabular}{cc}
    \includegraphics[width=0.3\textwidth]{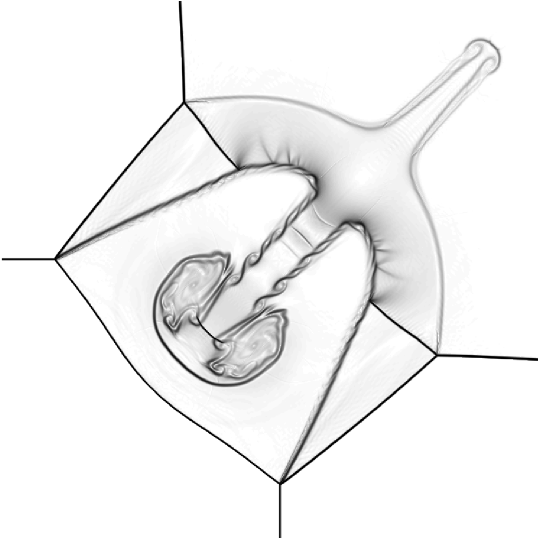} & \includegraphics[width=0.3\textwidth]{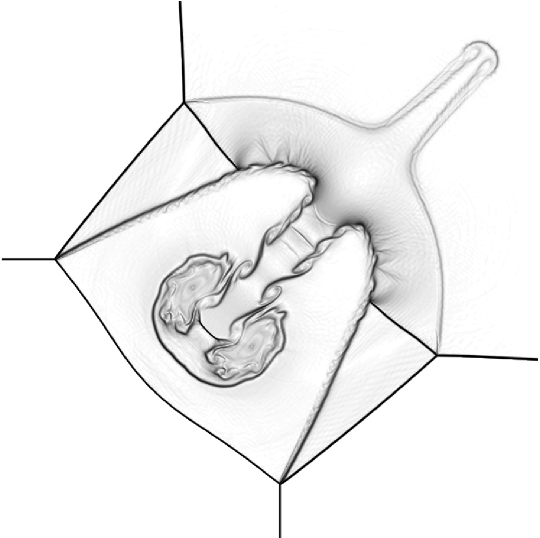} \\
    (a) WENO5, $2048\times2048$ & (b) CWENO5-LPR, $2048\times2048$ \\
    \includegraphics[width=0.3\textwidth]{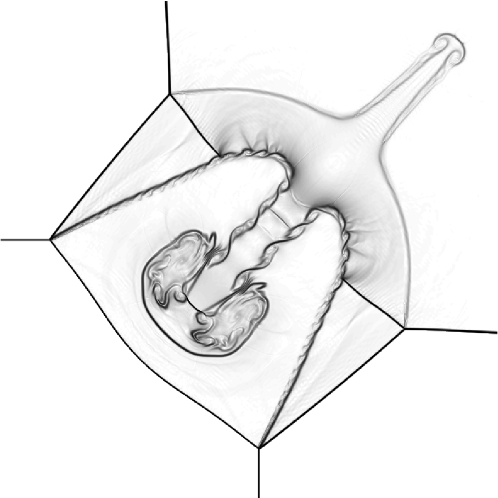} & \includegraphics[width=0.3\textwidth]{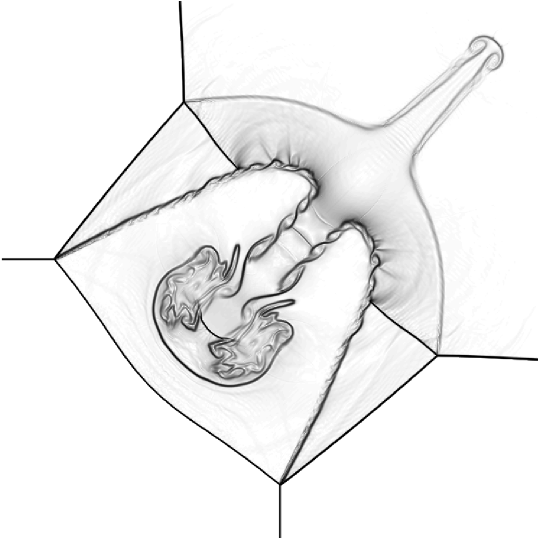} \\
    (c) CWENO5, $s_1=1$, $2048\times2048$ & (d) CWENO5, $s_1=2$, $2048\times2048$ \\
    \includegraphics[width=0.3\textwidth]{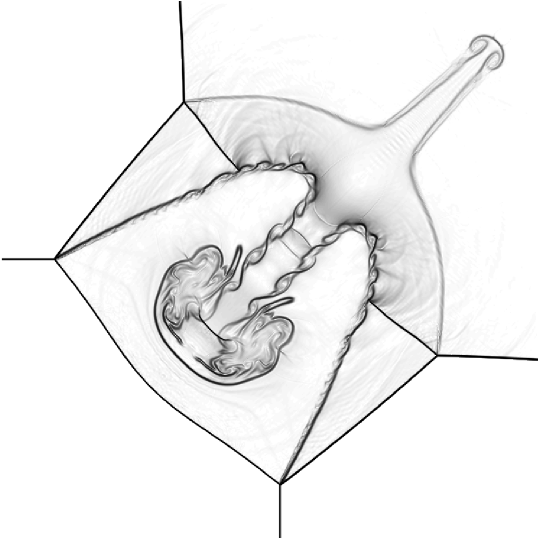} & \includegraphics[width=0.3\textwidth]{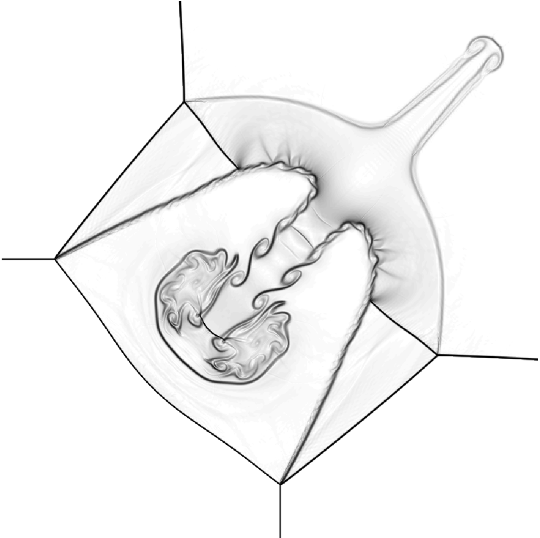}  \\
    (e) CWENO5, $s_1=3$, $2048\times2048$ & (f) WENO5, $2560\times2560$ \\
    \includegraphics[width=0.3\textwidth]{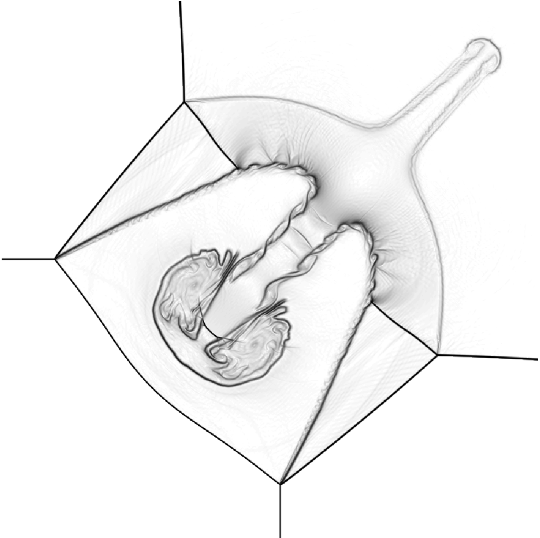} & \includegraphics[width=0.3\textwidth]{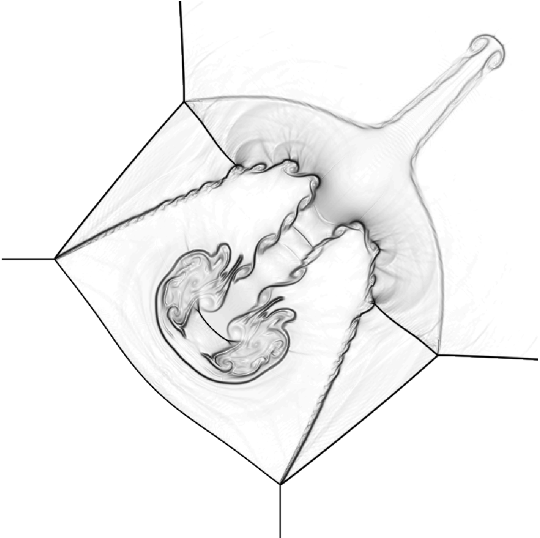}  \\
    (g) CWENO5-LPR, $2560\times2560$ & (h) CWENO5, $s_1=1$, $2560\times2560$ \\
    \includegraphics[width=0.3\textwidth]{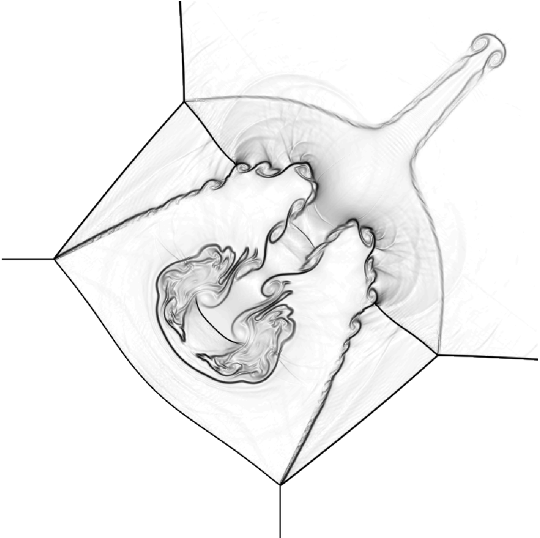} & \includegraphics[width=0.3\textwidth]{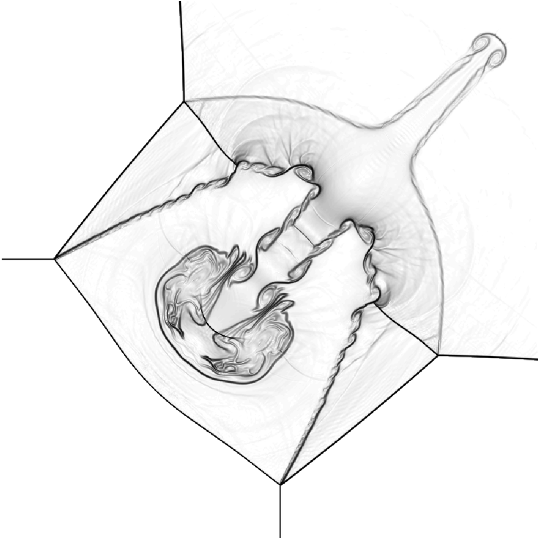}  \\
    (i) CWENO5, $s_1=2$, $2560\times2560$ & (j) CWENO5, $s_1=3$, $2560\times2560$
  \end{tabular}
  \caption{Example~7: enlarged view of turbulent zone (Schlieren plot).}
  \label{riemann}
\end{figure}

The simulation is performed taking $s_2=2$, with the final time $T=0.3$, $\textnormal{CFL}=0.4$, resolutions $2048\times2048$ and $2560\times2560$ and comparing the same schemes with the same parameters than in Example 6. The results are shown in Figure \ref{riemann}.

As shown in the results, in this case both CWENO5-LPR and CWENO5 have better resolution than WENO5, having in this case the former ones a similar resolution between them. The improvement of the resolution as $s_1$ increases is also observed in CWENO5 schemes. It is also important to point out that the results obtained with CWENO5 schemes in a resolution of $2048\times2048$ grid points are similar to those obtained with WENO5 and CWENO5-LPR in a resolution of $2560\times2560$ points.

\section{Conclusions}\label{cnc}

This paper proposes a novel WENO approach based on the computation of a global average weight as an additional measure of the global smoothness in a stencil. Such weight is then used to confer a stronger control involving the average between a reconstruction using the whole stencil, which ideally is performed when there is smoothness, and a reconstruction using properly half of the information from a smooth region, ideally performed when there is a discontinuity in the stencil.

The proposed scheme outperforms both the original WENO and CWENO-LPR schemes, which in the last case occurs in terms both of computational cost and error, as can be seen in our numerical experiments. We have also seen both theoretically and in practice that this approach overcomes whenever possible the well-known issues of the classical WENO schemes involving the loss of accuracy near smooth extrema and handling properly discontinuities, both quantitatively (theoretical results and numerical evidence) and qualitatively (best resolution near smooth extrema and shocks). Moreover, this approach can be used to easily tackle further issues regarding the WENO schemes, such as the presence of negative ideal weights for certain reconstructions (see for instance \cite{shu09} for further details involving reconstructions with negative weights) and the reconstruction at non-centered points.

Regarding the aforementioned considerations, we plan to use these schemes in several contexts, such as WENO reconstructions of derivatives and as a part of a more accurate boundary extrapolation algorithm, which will be illustrated in forthcoming works. We are also working with suitable weight designs for WENO schemes capable of attaining the optimal accuracy in presence of critical points, regardless of their order, planning also to extend such design and analysis to CWENO schemes. Finally, we are also considering to extend the idea behind the global average weight design to generalize our CWENO schemes in the context of unstructured grids.

\section*{Acknowledgments}
AB, PM and DZ are supported by Spanish MINECO project MTM2017-83942-P.  
RB is supported by Fondecyt project 1170473; CRHIAM, project CONICYT/FON\-DAP/15130015; and CONICYT/PIA/Concurso Apoyo a Centros Cient\'ificos y Tecnol\'{o}gicos 
 de Excelencia con Financiamiento Basal AFB170001. 
PM is also
supported by Conicyt (Chile),   pro\-ject PAI-MEC, folio 80150006.  
 DZ is also supported by Conicyt (Chile) through Fondecyt project 3170077.


\begin{thebibliography}{10}
  \providecommand{\url}[1]{{#1}}
  \providecommand{\urlprefix}{URL }
  \expandafter\ifx\csname urlstyle\endcsname\relax
  \providecommand{\doi}[1]{DOI~\discretionary{}{}{}#1}\else
  \providecommand{\doi}{DOI~\discretionary{}{}{}\begingroup
    \urlstyle{rm}\Url}\fi

\bibitem{SINUM2011}
  Aràndiga, F., Baeza, A., Belda, A.M., Mulet, P.: Analysis of {WENO} schemes
  for full and global accuracy.
  \newblock SIAM J.\ Numer.\ Anal. \textbf{49}(2), 893--915 (2011)

\bibitem{BaezaMuletZorio2016}
  Baeza, A., Mulet, P., Zorío, D.: High order weighted extrapolation for boundary conditions for finite difference methods on complex domains with Cartesian meshes.
  \newblock J.\ Sci.\ Comput. \textbf{69}(2), 170--200 (2016)

\bibitem{Capdeville2008}
  G. Capdeville: A central WENO scheme for solving hyperbolic conservation laws on non-uniform meshes.
  \newblock J. Comput. Physics. 227(5), 2977--3014 (2008)

\bibitem{CraveroSemplice}
  Cravero, I., Semplice, M.: On the accuracy of {WENO} and {CWENO} reconstructions of third order on nonuniform meshes.
  \newblock J.\ Sci.\ Comput. \textbf{67}, 1219--1246 (2015)

\bibitem{DonatMarquina96}
  Donat, R., Marquina, A.: Capturing shock reflections: An improved flux formula.
  \newblock J.\ Comput.\ Phys. \textbf{125}, 42--58 (1996)

\bibitem{Harten1987}
  Harten, A., Engquist, B., Osher, S., Chakravarthy, S.R.: Uniformly high order
  accurate essentially non-oscillatory schemes, {III}.
  \newblock J.\ Comput.\ Phys. \textbf{71}, 231--303 (1987)

\bibitem{Holoborodko}
  Holoborodko, P.: MPFR C++.
  \newblock http://www.holoborodko.com/pavel/mpfr/

\bibitem{JiangShu96}
  Jiang, G.S., Shu, C.W.: Efficient implementation of {Weighted} {ENO} schemes.
  \newblock J.\ Comput.\ Phys. \textbf{126}, 202--228 (1996)

\bibitem{KurganovTadmor}
  Kurganov, A., Tadmor, E.: Solution of two-dimensional {R}iemann problems for gas dynamics without Riemann problem solvers
  \newblock Numer.\ Methods  Partial Differential Equations  \textbf{18}, 584--608 (2002)

\bibitem{LevyPuppoRusso99}
Levy, D., Puppo, G., Russo, G.: Central {WENO} schemes for hyperbolic conservation laws.
\newblock ESAIM: Mathematical Modelling and Numerical Analysis \textbf{33}(3), 547--571 (1999)

\bibitem{LevyPuppoRusso}
  Levy, D., Puppo, G., Russo, G.: Compact {C}entral {WENO} schemes for multidimensional conservation laws.
  \newblock SIAM J. Sci. Comput. \textbf{22}(2), 656--672 (2000)

\bibitem{LiuOsherChan94}
  Liu, X-D., Osher, S., Chan, T.: Weighted essentially non-oscillatory schemes.
  \newblock J.\ Comput.\ Phys. \textbf{115}, 200--212 (1994)

\bibitem{MPFR}
  The GNU MPFR library.
  \newblock http://www.mpfr.org/

\bibitem{SchulzRinne}
  Schulz-Rinne, C. W.: Classification of the Riemann problem for two-dimensional gas dynamics.
  \newblock SIAM Journal on Mathematical Analysis \textbf{24}(1), 76--88 (1993)

\bibitem{shu98} 
Shu, C.-W.: Essentially non-oscillatory and weighted essentially non-oscillatory schemes for hyperbolic conservation laws. In 
Cockburn, B., Johnson, C., Shu, C.-W. and Tadmor, E.  (Quarteroni, A., Ed.):  Advanced Numerical Approximation of Nonlinear Hyperbolic Equations. 
  Lecture Notes in Mathematics, vol.~1697, Springer-Verlag, Berlin, 325--432 (1998) 

\bibitem{shu09} Shu, C.-W.:  High order weighted essentially nonoscillatory schemes 
 for convection dominated problems. SIAM Rev. {\bf 51}, 82--126 (2009) 

\bibitem{ShuOsher89}
  Shu, C.-W., Osher, S.: Efficient implementation of essentially non-oscillatory
  shock-capturing schemes.
  \newblock J.\ Comput.\ Phys. \textbf{77}, 439--471 (1988)

\bibitem{ShuOsher1989}
  Shu, C.-W., Osher, S.: Efficient implementation of essentially non-oscillatory
  shock-capturing schemes, {II}.
  \newblock J.\ Comput.\ Phys. \textbf{83}, 32--78 (1989)

\bibitem{zhangshu16} Zhang, Y.-T., Shu, C.-W.: ENO and WENO schemes. Chapter~5 in Abgrall, R.\ and Shu, C.-W.\ (eds.), 
   Handbook of Numerical Methods for Hyperbolic Problems Basic and Fundamental Issues. 
   Handbook of Numerical Analysis vol.~17, North Holland, 103--122 (2016) 
  
\end{thebibliography}
\end{document}